\RequirePackage{ifpdf}
\ifpdf 
\documentclass[pdftex]{sigma}
\else
\documentclass{sigma}
\fi

 \usepackage[mathscr]{euscript}
\input{xy}
\xyoption{all}

\newtheorem*{propn}{Proposition}
\newtheorem*{cartest}{Cartan's test}
\newtheorem*{cor}{Corollary}
\newtheorem*{lem}{Lemma}

 \theoremstyle{definition}
 \newtheorem*{cons}{Contact system}

\newtheorem*{defs}{Definitions}
\newtheorem*{defi}{Definition}
\newtheorem*{exam}{Example}
\newtheorem*{rems}{Remarks}
\newtheorem*{rem}{Remark}

\newcommand\uquiv{\begin{array}{c}=\\[-9.5pt]=\end{array}}
\newcommand\upt{{}^t}
\newcommand\bsl{\backslash}
\newcommand\checkd{\check{D}} 
\newcommand\checkh{\check{H}}
\newcommand\fbul{F^{\bullet}}
\newcommand{\upin}{\cup\hskip-4.95pt\raise1.5pt\hbox{$\scriptscriptstyle{|}$}}

\newcommand\wrvert{\lower1.25pt\hbox{$\wr$}{\scriptstyle\Vert}}
\newcommand\Vertwr{{\scriptstyle\Vert}\lower1.25pt\hbox{$\wr$}}

\newcommand\lrar[1]{\stackrel{#1}{\relbar\hskip-4pt\relbar\hskip-4pt\relbar\hskip-4pt\relbar\hskip-4pt\rightarrow}}
\newcommand\mrar[1]{\stackrel{#1}{\relbar\hskip-4pt\relbar\hskip-4pt\rightarrow}}

 \newcommand\srar[1]{\stackrel{#1}{\rightarrow}}
\newcommand\ol[1]{\bar{#1}}

\newcommand\ritem[1]{\item[{\rm #1}]}

\newcommand{\Gr}{{\rm Gr}}
\newcommand{\Hg}{{\rm Hg}}

 \newcommand\uh{{\mathbf{h}}}

\newcommand\Ad{\mathop{\rm Ad}\nolimits}
\newcommand\ad{\mathop{\rm ad}\nolimits}

 \newcommand\Aut{\mathop{\rm Aut}\nolimits}
\newcommand\End{\mathop{\rm End}\nolimits}
\newcommand\GL{\mathop{\rm GL}\nolimits}
 \newcommand\Grass{\mathop{\rm Grass}\nolimits}

\newcommand\rim{\mathop{\rm Im}\nolimits}

\newcommand\Sym{\mathop{\rm Sym}\nolimits}
\newcommand\rank{\mathop{\rm rank}\nolimits}
\newcommand\prim{{\rm prim}}

\newcommand{\limto}{\,{\displaystyle{\lim_{\raise2pt\hbox{$\scriptstyle\longrightarrow$}}}}}
\renewcommand{\part}{\partial}
\newcommand{\la}{{\lambda}}
\newcommand{\La}{{\Lambda}}
\newcommand{\Om}{{\Omega}}
\newcommand{\om}{{\omega}}

\newcommand\bb{\mathfrak{b}}
\newcommand\scri{\mathfrak{i}}
\newcommand\pp{\mathfrak{p}}
\newcommand\fv{\mathfrak{v}}
\newcommand\fw{\mathfrak{w}}
\newcommand\fadv{\hbox{ad-}\fv}

\newcommand{\scrA}{{\frak{a}}}

\newcommand{\scrF}{{\mathscr{F}}}
\newcommand{\scrG}{{\mathscr{G}}}

\newcommand{\scrI}{{\mathscr{I}}}

\newcommand{\scrM}{{\mathscr{M}}}

\newcommand{\scrO}{{\mathscr{O}}}

\newcommand{\scrU}{{\mathscr{U}}}

\renewcommand{\P}{\mathbb{P}}
\newcommand{\Q}{\mathbb{Q}}

\newcommand{\R}{\mathbb{R}}

\newcommand{\Z}{\mathbb{Z}}
\newcommand{\C}{\mathbb{C}}

\newcommand{\pref}[1]{(\ref{#1})}

\newcommand{\lab}{\label}
\newcommand{\lra}[1]{\left\langle#1\right\rangle}
\newcommand{\lrc}[1]{\left\{ #1\right\}}
\newcommand{\lrb}[1]{\left[ #1\right]}
\newcommand{\lrp}[1]{\left(#1\right)}
\newcommand{\Hom}{\mathop{\rm Hom}\nolimits}

\newcommand{\opplus}{\mathop{\oplus}\limits}
\newcommand{\ottimes}{\mathop{\otimes}\limits}

\newcommand{\codim}{\mathop{\rm codim}\nolimits}

\newcommand{\hensp}[1]{\enspace\hbox{#1}\enspace}

\numberwithin{equation}{section}

\begin{document}

\allowdisplaybreaks

\renewcommand{\thefootnote}{$\star$}

\renewcommand{\PaperNumber}{087}

\FirstPageHeading

\ShortArticleName{Variations of Hodge Structure}

\ArticleName{Variations of Hodge Structure Considered\\ as an Exterior Dif\/ferential
System:\\ Old and New Results\footnote{This paper is a
contribution to the Special Issue ``\'Elie Cartan and Dif\/ferential Geometry''. The
full collection is available at
\href{http://www.emis.de/journals/SIGMA/Cartan.html}{http://www.emis.de/journals/SIGMA/Cartan.html}}}

\Author{James CARLSON~$^\dag$, Mark GREEN~$^\ddag$ and Phillip GRIFFITHS~$^\S$}

\AuthorNameForHeading{J.~Carlson, M.~Green and  P.~Grif\/f\/iths}

 \Address{$^\dag$~Clay    Mathematics Institute, United States}
 \EmailD{\href{mailto:jxxcarlson@mac.com}{jxxcarlson@mac.com}}

 \Address{$^\ddag$~University of California, Los Angeles, CA, United States}
 \EmailD{\href{mailto:mlgucla@gmail.com}{mlgucla@gmail.com}}

\Address{$^\S$~The Institute for Advanced Study, Princeton, NJ, United States}
\EmailD{\href{mailto:pg@ias.edu}{pg@ias.edu}}

\ArticleDates{Received April 20, 2009, in f\/inal form August 31, 2009;  Published online September 11, 2009}

 \newcounter{demo}[equation]
   \renewcommand{\thedemo}{\thesection.\arabic{equation}}

 \newenvironment{demo}{
 \addtocounter{equation}{1}\refstepcounter{demo} \setcounter{demo}{\value{equation} }
{\medbreak\noindent (\thesection.\arabic{demo})\enspace   }}{\medbreak}

\Abstract{This paper is a survey of the subject of variations of Hodge
structure (VHS)
considered as exterior dif\/ferential systems (EDS).  We review
developments
over the last twenty-six years, with an emphasis on some key examples.
In the penultimate section we present some new results on the
characteristic cohomology
of a homogeneous Pfaf\/f\/ian system.  In the last section we discuss how
the integrability
conditions of an EDS af\/fect the expected dimension of an integral
submanifold.
The paper ends with some speculation on EDS and
Hodge conjecture for Calabi--Yau manifolds.}


\Classification{14C30; 58A15}

 \begin{flushright}
\begin{minipage}{15cm}
\small\it
A portion of this paper was presented by the third named author during the Conference on Exterior Differential Systems and Control Theory held at the Mathematical Science Research Institute in Berkeley.  This conference was held in honor of Robby Gardner, whose contributions to both exterior differential systems and control theory were of the greatest significance.  The authors would like to thank the organizers for putting together
        the conference and would like to dedicate this paper to the memory of
Robby Gardner.
\end{minipage}
\end{flushright}

 {
 \tableofcontents }

%

  \renewcommand{\thefootnote}{\arabic{footnote}}
\setcounter{footnote}{0}

  \section{Introduction}

  Hodge theory provides the basic invariants of a complex algebraic variety.  The two central open problems in the subject, the conjectures of Hodge and of Bloch--Beilinson, relate Hodge theory to the geometry/arithmetic of a complex algebraic variety.

  The space of all polarized Hodge structures of weight $n$ and with a given sequence $\uh = (h^{n,0},h^{n-1,1},\dots,h^{0,n})$, $h^{p,q}=h^{q,p}$, of Hodge numbers forms naturally a homogeneous complex manifold $D_{\uh}$, and the moduli space $\scrM_{\uh}$ of equivalence classes of polarized Hodge structures is a~quotient of $D_{\uh}$ by an arithmetic group acting properly
   discontinuously.   In what we shall call the {\it classical case} when the weights $n=1$ or $n=2$ and $h^{2,0}=1$, $D_{\uh}$  is a bounded symmetric domain and $\scrM_\uh$ is a quasi-projective variety def\/ined over
    a number f\/ield\footnote{Henceforth we shall drop reference to the $\uh$ on $D$ unless it is needed.}.  In this case the relation between the  Hodge theory and geometry/arithmetic of a variety is an extensively developed and deep subject.

   In the non-classical, or what we shall refer to as the {\it higher weight}, case\footnote{It being understood that when $n=2$ we have $h^{2,0}\geqq 2$.} the subject is relatively less advanced.  The fundamental dif\/ference between the classical and higher weight cases is that in the latter case the Hodge structures associated to a family of algebraic varieties satisfy a universal system $I\subset T^\ast D$ of dif\/ferential equations.  In this partly expository paper we will discuss the system $I$ from the perspective of {\it exterior differential systems} (EDS's) with the three general  objectives:
      \begin{itemize}\itemsep=0pt
   \ritem{$(i)$} To summarize some of  what is known about $I$ from an EDS perspective.
   \ritem{$(ii)$} To def\/ine and  discuss the ``universal characteristic cohomology" associated to a homogeneous Pfaf\/f\/ian system in the special case of a variation of Hodge structures.

    \ritem{$(iii)$} To discuss and illustrate the   question ``How must expected dimension counts be modif\/ied for integral manifolds of the system~$I$?"
   \end{itemize}

One overarching objective is this: When one seeks to extend much of the rich classical theory, including the arithmetic aspects and the connections with automorphic forms, the various compactif\/ications  of $\scrM$ and the resulting boundary cohomology, the theory of Shimura varieties, etc., the fact that families of Hodge structures arising from geometry are subject to dif\/ferential constraints seems to present {\it the} major barrier.  Perhaps by better understanding the structure of these dif\/ferential constraints, some insight might be gained   on how at least some aspects of the classical theory might be extended.  We are especially interested in properties of variations of Hodge structure that are {\it not} present in the classical case, as these may help to indicate what needs to be better understood to be able to extend the classical case to higher weight Hodge structures.

   In more detail, in Section \ref{2a} we will review the def\/initions and establish notations for polarized Hodge structures, period domains and their duals, and the inf\/initesimal period relation, which is the basic exterior dif\/ferential system studied in this paper.  In Section~\ref{2.b} we recall some of the basic def\/initions and concepts from the theory of exterior dif\/ferential systems.

   In Section \ref{3} we discuss the basic exterior dif\/ferential system whose integral manifolds def\/ine variations of Hodge structure.   The basic general observation is that the integral elements are given by {\it abelian}  subalgebras
\[
   \scrA\subset \scrG^{-1,1}\subset \scrG_\C
   \]
   of the complexif\/ied Lie algebras of the symmetry group of the period domain.  We then go on in Section~\ref{3.a} to discuss in some detail two important examples where the EDS given by the inf\/initesimal period relation may be integrated  by elementary methods; both of these have been discussed in the literature and here we shall summarize, in the context of this paper  those results.  Then in Section~\ref{3.b}  we shall discuss an example,  the f\/irst we are aware of in the literature, where the Cartan--K\"ahler theoretic  aspects of exterior dif\/ferential systems are applied to the particular EDS arising from the inf\/initesimal period relation.  Finally, in Section~\ref{3.c} we study the derived f\/lag of the inf\/initesimal period relation $I$.  One result is that if all the Hodge numbers $h^{p,q}$ are non-zero, then the derived f\/lag of $I$ terminates in zero and has length $m$ where $m\leqq \frac{\log n}{\log 2}$, $n$~being the weight of the Hodge structure.

   Over the years there have been a number of studies of the EDS given by the inf\/initesimal period relation.
   Here we mention \cite{A,C, CD, CKT, CS, CT, CT2} and \cite{M}.  Section \ref{3} should be considered as an introduction to those works.  In particular the   paper~\cite{M}, which builds on and extends the earlier works, contains a def\/initive account of the bounds on the dimension and rigidity properties of maximal integral elements.  At the end of Section \ref{3.c} we shall comment on some interesting questions that arise from~\cite{M} and~\cite{A}  as well as providing a brief guide to the earlier works referred to above.

   In Section~\ref{4}, we f\/irst discuss some general aspects of homogeneous Pfaf\/f\/ian systems, including expressing the invariant part of their {\it characteristic cohomology} in terms of a Lie algebra cohomology construction\footnote{Here, we recall that the characteristic cohomology of an exterior dif\/ferential system is the de~Rham cohomological construction  that leads to cohomology groups that induce ordinary de~Rham cohomology on integral manifolds of the EDS.  The precise def\/inition is recalled in Section~\ref{4} below.}.  We then turn to the group invariant characteristic cohomology of period domains.  Here there is a very nice question
   \begin{demo}\lab{1.1}
 {\it Is the invariant part of the characteristic  cohmology of a period domain   generated by the Chern forms of the Hodge bundles?}
   \end{demo}

   In the classical case the answer is positive and may be deduced from what is known in the literature; we will carry this out below.

   In the non-classical case when the Pfaf\/f\/ian system $I$ associated to a variation of Hodge structures is non-trivial, new and interesting issues arise.  It seems likely that the question will have an af\/f\/irmative answer; this will be the topic of a separate
   work\footnote{This question has now been answered in the af\/f\/irmative and the proof will appear in a separate publication.}.

   In this paper we will establish two related results.  The f\/irst is that we will show that the invariant forms modulo the {\it algebraic} ideal generated by $I$ are all of type $(p,p)$.  A consequence is that on the complex of invariant forms the Lie algebra cohomology dif\/ferential $\delta=0$.  This is analogous to what happens in the Hermitian symmetric case.  However, in the non-classical case there are always more invariant forms than those generated by the Chern forms of the Hodge bundles, and the {\it integrability conditions}; i.e.\ the full {\it differential} ideal generated by $I$ must be taken into account.  This involves subtle issues in representation theory and, as mentioned above, this story will be reported on separately.

   What we will prove here is that the integrability conditions imply topological conditions in the form of new relations among the Chern classes of a VHS.  Denoting by $\scrF^p$ the Hodge f\/iltration bundles we show that the Chern forms satisfy
   \[
   c_i(\scrF^p)c_j(\scrF^{n-p})=0\qquad \hensp{if} \quad i+j > h^{p,n-p}.\]
   What we see then is that this new algebro-geometric information has resulted from EDS considerations.

      The question \pref{1.1} would have the following algebro-geometric consequence:
      First, for any discrete group $\Gamma$ acting properly discontinuously on a period domain $D$, the invariant characteristic cohomology $H^\ast_{\scrI}(D)^{G_\R}$ induces characteristic cohomology on the quotient $\Gamma\bsl D$.  A global variation of Hodge structure is given by a period mapping
    \begin{equation}\lab{1.2}
   f:  \ S\to \Gamma\bsl D,
   \end{equation}
    where  $S$ is a smooth, quasi-projective algebraic variety.  Since \pref{1.2} is an integral manifold of the canonical    EDS on $\Gamma\bsl D$, the invariant characteristic cohomology induces ordinary cohomology
    \[
   f^\ast\lrp{H^\ast_\scrI(D)^{G_\R}}\subset H^\ast(S) .\]
    We may think of $H^\ast_\scrI(D)^{G_\R}$ as {\it universal characteristic cohomology} for the exterior dif\/ferential system corresponding to the inf\/initesimal period relation, in that it induces ordinary cohomology on the parameter space for variations of Hodge structures  {\it irrespective of the particular group} $\Gamma$. We may thus call it {\it the} universal characteristic cohomology.  A positive answer to the question would f\/irst of all  imply that the universal characteristic cohomology is generated by the Chern classes of the Hodge bundles over $S$.  Although, perhaps not surprising to an algebraic geometer  this would be  a satisfying result.

    In Section~\ref{5} we turn to the interesting question
    \begin{demo}\lab{1.3}
    {\it How must one correct expected dimension counts in the presence of differential constraints?}
    \end{demo}

      Specif\/ically, given a manifold $A$ and submanifold $B\subset A$,  for a ``general" submanifold $X\subset A$ where $\dim X + \dim B\geqq \dim A$, we  will have
\begin{equation}\lab{1.4}
 \codim_B(X\cap B)=\codim_A B.
 \end{equation}
Thus, the RHS of this equation may be thought of as the ``expected codimension" of  $X\cap B$ in~$B$.  If $X\cap B$ is non-empty, the actual codimension is no more than  the expected codimension.

Suppose now that there is a distribution $W\subset TA$ and $X$ is constrained to have $TX\subset W$.   Then how does this af\/fect the expected dimension counts?  In case $W$ meets $TB$ transversely, one sees immediately that the ``expected codimension" counts decrease.  Taking integrability into account gives a further correction.  Rather than trying to develop the general theory, in Section~\ref{5} we shall discuss one particularly interesting special case.

This case concerns {\it Noether--Lefschetz loci}.  Here we denote by $W_I\subset TD$ the  distribution~$I^\bot$; integral manifolds of $I$ have their tangent spaces lying in~$W_I$.  Given $\zeta\in H_\R$, there is  a homogeneous sub-period-domain $D_\zeta\subset D$, def\/ined as the set of polarized Hodge structures of weight $n=2m$ where $\zeta\in H^{m,m}$.  We have
\[
\codim_D D_\zeta=h^{(2m,0)}+\dots+ h^{(m+1,m-1)},\]
and the distribution $W_I$ on $D$ meets $TD_\zeta$ transversely.  For a variation of Hodge structure, given by an integral manifold
\[ f: \ S\to D\]
of the canonical system $I$ on $D$, the Noether--Lefschetz locus $S_\zeta\subset S$ is given by
\[
S_\zeta=f^{-1}(D_\zeta).\]
In algebro-geometric questions, $\zeta$ is usually taken to be a rational vector, but that will not concern us here.   The ref\/ined codimension estimate given by $W_I$ alone, i.e., without taking integrability conditions into account, is
\[
\codim_S S_\zeta\leqq h^{m-1,m+1}.\]
In the case $m=1$, which algebro-geometrically ref\/lects studying codimension one algebraic cycles,  the distribution $W_I$ does not enter and the estimate is classical, especially in the study of curves on an algebraic surface.  In the case $m\geqq 2$, it is non-classical and seems only recently to have been discussed  in the literature (cf.~\cite{GG2,V, O1,O2}).
    When the integrability conditions are taken into account, the above codimension estimate is ref\/ined to
    \[
    \codim_S S_\zeta\leqq h^{m-1,m+1}-\sigma_\zeta,\]
    where $\sigma_\zeta$ is a non-negative quantity constructed from $\zeta$ and the integral element of $I$ given by the tangent space to $f(S)$ at the point in question.  Assuming the Hodge conjecture, the above would say that ``there are more algebraic cycles than a na\"ive dimension count would suggest".

    In addition to establishing the above inequality, we will show that it is an equality in a~signif\/icant example, namely, that given by a   hypersurface $X\subset \P^5$ of degree $d\geqq 6$ and which is general among those containing a $2$-plane. This  indicates that there is no further general estimate.

    We conclude this section by analyzing the case of Calabi--Yau fourfolds, where the quantity $h^{3,1}-\sigma_\zeta$ has a particularly nice   interpretation, including an interesting arithmetic consequence of the Hodge conjecture.

   \section{Preliminaries}

   \subsection[Period domains]{Period domains\footnote{The general reference for this section is the book \cite{C-MS-P}.}}\label{2a}

   Let $H$ be a $\Q$-vector space.  A {\it Hodge structure} (HS) {\it of weight} $n$ is given by any of the following equivalent data:

   \begin{demo}   \label{2.1} A {\it Hodge decomposition}
   \[
   \left\{\begin{array}{l}
   H_\C = \opplus_{p+q=n} H^{p,q},\\[11pt]
   H^{q,p} =\ol{H}^{p,q}.\end{array}\right.
   \]
   \end{demo}

   \begin{demo} \lab{2.2} A {\it Hodge filtration} $F^n\subset F^{n-1}\subset\dots \subset F^0= H_\C$ where for each~$p$
   \[
   F^p\oplus \ol{F}^{n-p+1}\stackrel{\sim}{\to} H_\C  .\]
   \end{demo}
   \begin{demo} \lab{2.3}
   A homomorphism of real Lie groups
   \[
   \varphi: \  \C^\ast\to \GL(H_\R)\]
   of weight $n$ in the sense that for $z\in \C^\ast$, $\la\in \R^\ast$
   \[\varphi(\la z)=\la^n \varphi(z).\]
   \end{demo}
   \medbreak\noindent
   The relation between \pref{2.1} and \pref{2.2} is
   \[
   \left\{
   \begin{array}{l}
   F^p   = \opplus_{p'\geqq p} H^{p',n-p'},\\[11pt]
   H^{p,q}   = F^p\cap \ol{F}^q .\end{array}\right.\]
   The relation between \pref{2.1} and \pref{2.3} is
   \[
 \varphi(z)u=z^p\ol{z}^q u ,\qquad  u\in H^{p,q} .\]
   This means: The element $\varphi(z)\in \GL(H_\C)$ just  given    lies in the subgroup $\GL(H_\R)$ of $\GL(H_\C)$ and has weight $n$.

   If we restrict $\varphi$ to the maximal compact subgroup $S^1=\{z\in\C^\ast:|z|=1\}$ of $\C^\ast$, then for $z\in S^1$, $u\in H^{p,q}$
   \[
   \varphi(z)u = z^{p-q}u ,\]
    and this shows how to recover the Hodge decomposition as the $z^{p-q}$ eigenspace of $\varphi$ restricted to $S^1$.
    The {\it Weil operator} is def\/ined by
    \[
    C=\varphi(\sqrt{-1}) .\]

    We def\/ine the Hodge numbers $h^{p,q}:=\dim H^{p,q}$, and we set $f^p :=\sum_{p^1\geqq p} h^{p',n-p'}$.

    Since the Hodge f\/iltration point of view will be the dominant one in this paper, we shall denote a Hodge structure by $(H,F)$.

    Now let
    \[
    Q: \ H\otimes H\to \Q\]
    be a non-degenerate bilinear form satisfying $Q(u,v)=(-1)^nQ(v,u)$ for $u,v\in H$.  A Hodge structure $(H, F)$ is {\it polarized} by $Q$ if the Hodge--Riemann bilinear relations
    \begin{equation} \lab{2.4}
    \left\{
    \begin{array}{l}
    Q(F^p,F^{n-p+1})=0,\\[4pt]
   Q(Cu,\ol{u}) >0 \quad \hensp{for} \ \ 0\ne u\in H_\C \end{array}\right.
   \end{equation}
   are satisf\/ied.
   The f\/irst of these is equivalent to $F^{n-p+1}=(F^p)^\bot$.  A polarized Hodge structure (PHS) will be denoted by $(H,Q,F)$.

   In the def\/inition \pref{2.3}, for a polarized Hodge structure we need to restrict $\varphi$ to $S^1$ in order to preserve, and not just scale, the polarization.

   \begin{defs}\rm  $(i)$ A {\it period domain} $D$ is given by the set of polarized Hodge structures $(H,Q,F)$ with given Hodge numbers $h^{p,q}$.  $(ii)$ The {\it compact dual} $\check{D}$ is given by all f\/iltrations $F$ with $\dim F^p=f^p$ and which satisfy the f\/irst bilinear relation $F^{n-p+1}=(F^p)^\bot$ in~\pref{2.4}.
   \end{defs}

   We shall denote by $G$ the $\Q$-algebraic group $\Aut(H,Q)$, and by $G_\R$ and $G_\C$ the corresponding real and complex forms.  It  is elementary that $G_\R$ acts transitively on $D$, and choosing a reference point $F_0\in D$ we have
   \[
   D\cong G_\R/V,
   \]
   where $V$ is the compact subgroup of $G_\R$ preserving the Hodge decomposition $H_\C=\opplus_{p+q=n}H^{p,q}_0$ corresponding to $F_0$.  In terms of \pref{2.3} we note that
   \[
   \varphi(S^1)\subset V ;
   \]
   in fact, $V$ is the centralizer of the circle $\varphi(S^1)$ in $G_\R$.

     The complex Lie group $G_\C$ acts transitively on the compact dual $\check{D}$, and choosing a reference point $F_0 \in D$ as above we have
   \[
   \check{D}\cong G_\C/B,
   \]
   where $B\subset G_\C$ is a parabolic subgroup.  We have
   \[
   V=G_\R\cap B .\]
   The compact dual is a projective algebraic variety def\/ined over $\Q$.  In fact we have an obvious inclusion
   \begin{equation}\lab{2.5}
   \check{D}\subset \prod^{\left[\frac{n+1}{2}\right]}_{p=n} \mathrm{Grass}(f^p,H_\C)
   \end{equation}
   and we may embed $\check{D}$ in a $\P^N$ by means of the Pl\"ucker coordinates of the f\/lag subspaces $F^p\subset H_\C$.

   Hodge structures and polarized Hodge structures are functorial with respect to the standard operations in linear algebra.  In particular, a Hodge structure $(H,F)$ induces a Hodge structure of weight zero on  $\End(H)$ where
   \begin{equation}\lab{2.6}
   \End(H_\C)^{r,-r}= \left\{ A\in \End(H_\C):A(H^{p,q})\subset H^{p+r,q-r}\right\} .
   \end{equation}
   A polarized Hodge structure $(H,Q,F_0)$ induces a polarized Hodge structure on the Lie algebra~$\scrG$ of~$G$, where $\scrG^{r,-r}$ is given by \pref{2.6} and the polarization is induced by the Cartan--Killing form.  For later use, we note that from~\pref{2.6} we have
   \begin{equation}\lab{2.7}
   \lrb{\scrG^{r ,-r},\scrG^{s,-s}}\subseteq \scrG^{r+s,-(r+s)} .\end{equation}

   If we recall the natural identif\/ication
   \[
   T_{F^p} \Grass(f^p,H_\C)\cong \Hom(F^p,H_\C/F^p)
   \]
   it follows that the Lie algebra of $B$ is
   \[
\bb = \opplus_{r\geqq 0}\scrG^{r,-r} .\]
   The subalgebra
   \[
   \pp = \opplus_{r>0}\scrG^{-r,r}\]
   gives a  complement to $\bb$ in $\scrG_\C$ leading to the natural identif\/ication
   \begin{equation}\lab{2.8}
   T_{\fbul} \check{D} \cong \pp .
   \end{equation}
   By \pref{2.7} the subspace
   \[
   \scrG^{-1,1}\subset \pp\]
   is $\Ad B$-invariant
   and therefore def\/ines a $G_\C$-invariant distribution
   \[
   W_I\subset T\check{D},
   \]
   and, by orthogonality, a Pfaf\/f\/ian system
   \[
   I\subset T^\ast \check{D} .\]
   The sub-bundle $I$ restricts to a $G_\R$-invariant sub-bundle $I\subset T^\ast D$.

   \begin{defi}  The Pfaf\/f\/ian system $I$ is called the {\it infinitesimal period relation}.\end{defi}

   It is this exterior dif\/ferential system that we shall discuss in this paper.

   In the literature the distribution $W_I\subset  TD$ is frequently referred to as the {\it horizontal sub-bundle}.

   \subsection[Exterior dif\/ferential systems (EDS)]{Exterior dif\/ferential systems (EDS)\footnote{General references  for this section are   books \cite{IL} and  \cite{BCGGG}; especially the former contains essentially all the background needed for this work.}}\lab{2.b}

    Although the subject is usually discussed in the smooth category, here we shall work complex-analytically.  A {\it Pfaffian system} is given by a holomorphic sub-bundle
      \[
      I\subset T^\ast M\]
      of the cotangent bundle of a complex manifold.   Associated to $I$ is the {\it differential ideal}
      \[
      \scrI\subset \Om^\bullet_M\]
      generated by the holomorphic sections of $I$ together with their exterior derivatives.   We shall assume that the values of the sections of $\scrI$ generate a sub-bundle of $\La^\bullet T^\ast M$; i.e., $\scrI$ is locally free.    An {\it integral manifold}, or just  an ``integral", of $I$ is given by a complex manifold $N$ and a~holomorphic immersion
      \[
      f: \ N\to M\]
      such that
      \begin{equation} \lab{2.9}
      f^\ast(\scrI)=0 .\end{equation}
      If we denote by
      \[
      W_I = I^\bot\subset TM\]
      the holomorphic distribution associated to $I$, the condition \pref{2.9} is equivalent to
      \[
      f_\ast(TN)\subset W_I .\]

      An important invariant associated to a Pfaf\/f\/ian system is its {\it derived flag}.  The exterior derivative induces a bundle map
      \[
      \delta : \ I\to \La^2T^\ast M/ I\wedge T^\ast M    , \]
      and recalling our assumption that $\delta$ has constant rank we set
      \[
      I_{[1]} = \ker \delta .\]
      This is again a Pfaf\/f\/ian system, and continuing in this way leads to the derived f\/lag in $T^\ast M$
      \begin{equation}\lab{2.10}
      I\supset I_{[1]} \supset I_{[2]} \supset \dots \supset I_{[m]} =I_{[m+1]} = \dots = I_{[\infty]} .
      \end{equation}
      Here, $I_{[\infty]}$ is the largest integrable or Frobenius subsystem of $I$.

      Dually, for the distribution we denote by
      \[
      W^{[1]}_I=W_I + [W_I,W_I]\]
      the distribution generated by $W_I$ and the brackets of sections of $W_I$.  Continuing in this way we obtain the
   f\/lag in $TM$ dual to \pref{2.10}
   \[
   W_I\subset W^{[1]}_I\subset\dots \subset W_I^{[m]} = W_I^{[m+1]} = \dots = W_I^{[\infty]} .
   \]
   We say that $I$ is {\it bracket generating} in case $W^{[\infty]}_I=TM$, or equivalently $I^{[\infty]}=(0)$.  In this case, by the holomorphic version of Chow's theorem we may connect any two points of $M$ by a chain of holomorphic discs that are integral curves of $W_I$.

   Two central aspects of the theory of exterior dif\/ferential systems are $(i)$ {\it regular and ordinary  integral elements} and the {\it Cartan--K\"ahler theorem}, and $(ii)$ {\it prolongation} and {\it involution}.  For the f\/irst, an {\it integral element} for $\scrI$ is given by a linear subspace $E\subset T_x M$ such that
   \[
   \theta_E =: \theta\mid_E=0\]
   for all $\theta\in \scrI$.  We may think of $E$ as an inf\/initesimal solution of the EDS.  Denote by
   \[
   \pi: \ G_p(TM)\to M\]
   the bundle whose f\/ibre $\pi^{-1}(x)=\Gr_p(T_x M)$ over $x\in M$ is the Grassmannian of $p$-planes in~$T_x M$.  In $G_p(TM)$ there is the complex analytic subvariety
   \[G_p (I) \subset G_p(TM)\]
   of integral elements def\/ined by the Pfaf\/f\/ian system $I$.

   Integral elements are constructed one step at a time by solving linear equations.  For an integral element $E\in G_p(I)$, we def\/ine the {\it polar space}
   \[
   H(E) = \lrc{
   v\in T_xM:\hbox{span}\{v,E\}\hbox{ is an integral element}} .\]
   The equations that def\/ine $H(E)$
   \[
   \lra{ \theta(x),v\wedge E}=0\qquad \hbox{for all} \quad \theta\in \scrI^{p+1}\]
   are linear in $v$, and we measure their rank by def\/ining
   \begin{gather*}
   r(E)=\codim H(E)
 =  \dim \P(H(E)/E) .
   \end{gather*}

   Given $E_0\in G_p(I)$ we choose a $p$-form $\Om$ such that $\Om_{E_0}\ne 0$.  For $\theta\in \scrI^p$ and $E\in G_p(TM)$ near $E_0$, for each $\varphi\in \scrI^p$ we write
   \[
   \theta_E=f_\theta(E)\Om_E .\]
   Then  $G_p(I)$ is locally def\/ined by the analytic equations $f_\theta(E)=0$; we say that $E_0$ is {\it regular} if~$G_p(I)$ is smoothly def\/ined by these equations.

   Given $E\in G_p(I)$ we choose a generic f\/lag
   \[
   0\subset E_1\subset\dots\subset E_{p-1}\subset E \]
   and def\/ine $c_i$ to be the rank of the polar equations of $E_i$.
   A central result is
   \begin{cartest}\qquad

   $(i)$ We have
   \[
   \mathop{\codim}  G_p(I)\geqq c_1+\dots + c_{p-1} .\]

   $(ii)$ If equality holds, then for $k<p$ each $E_k$ is regular.
    \end{cartest}

   By def\/inition we say that $E$ is {\it ordinary} if equality holds in Cartan's test.  We say that $I$ is {\it involutive} near an ordinary integral element.

   The {\it Cartan--K\"ahler existence theorem} states that an ordinary integral element is locally tangent to an integral manifold.  It goes further to say ``how many" local integral manifolds there are.  This will be brief\/ly mentioned below.

   If $G_p(I)$ is a submanifold near $E$, but $E$ is not ordinary one needs to {\it prolong} the EDS to be able to use the Cartan--K\"ahler theorem in order to construct local integral manifolds.  To explain this, we f\/irst observe that there is a canonical Pfaf\/f\/ian system
   \[
   J\subset T^\ast G_p(TM)
   \]
   def\/ined by writing points of $G_p(TM)$ as $(x,E)$ and setting
   \[
   J_{(x,E)} = \pi^\ast (E^\bot) .\]
   For any immersion $f:N\to M$ where $\dim N=p$, there is a canonical lift
   \[
   \xymatrix{
   & G_p(TM) \ar[d]^{\pi}\\
   N\ar[ur]^{ f^{(1)}} \ar[r]^f & M}
   \]
   where  for $y\in N$ we set  $f^{(1)} (y) = (f(y),f_\ast T_yN)$.  This lift is an integral manifold of $J$.  The converse is true provided the $p$-dimensional integral manifold of $J$ projects to an immersed $p$-dimensional submanifold
   in~$M$.

   We def\/ine the $1^{\rm st}$ {\it prolongation} $(M^{(1)}_I,I^{(1)})$ of $(M,I)$ by taking $M_I^{(1)}$ to be the smooth locus of $G_p(I)$ and def\/ining $I^{(1)}$ to be the restriction of $J$ to $M^{(1)}_I$.  An integral manifold of $(M,I)$ gives one of $(M^{(1)}_I,I^{(1)})$; the converse holds in the sense explained above.  The {\it Cartan--Kuranishi theorem} states that, with some technical assumptions, after a f\/inite number of prolongations an  exterior dif\/ferential system either becomes involutive or else empty (no solutions).

   From the above we may infer the following
     \begin{demo}
   \lab{2.11}
   {\it Let $f:N\to M$ be an integral manifold of a Pfaffian system.  If, at a general point $y\in N$, $f_\ast(T_yN)$ is not an ordinary integral element then one of the following holds:
   \begin{itemize}\itemsep=0pt
   \ritem{(A)} $f_\ast (TN)$ lies in the singular locus of $G_p(I)$; or
   \ritem{(B)} the prolongation $f^{(1)}:N\to M_I^{(1)}$ is an integral manifold of $I^{(1)}$.
   \end{itemize}
   }
   \end{demo}

 We may iterate this for the successive prolongations.  Since the failure to be involutive means that $f:N\to M$ satisf\/ies additional dif\/ferential equations not in $\scrI$, we have the following
   \begin{demo}
   \lab{2.12}
   {\bf Conclusion.}  {\it If the assumption of $(\ref{2.11})$ holds, then $f\!:\!N\!\to\! M$ satisfies additional differential equations that are canonically associated to $I$ but are not in the differential ideal $\scrI$.}
   \end{demo}

   \section[The exterior dif\/ferential system associated to a variation of Hodge structure]{The exterior dif\/ferential system associated\\ to a variation of Hodge structure}\lab{3}

   In this section, unless mentioned otherwise we shall work locally.  Let $D$ be a period domain and $I\subset T^\ast D$ the inf\/initesimal period relation with corresponding horizontal distribution $W_I\subset TD$.  We recall that both $I$ and $W_I$ are the restrictions to $D\subset \checkd$ of similar structures over the compact dual of $D$.

   \subsection{Elementary examples}\lab{3.a}

   \begin{defi} A {\it variation of Hodge structure} (VHS) is given by an integral manifold $f:S\to D$ of $I$.
  \end{defi}

  In the classical case where $D$ is a bounded symmetric domain, it is well known that holomorphic mappings from quasi-projective varieties to quotients $\Gamma\bsl D$ of $D$ by arithmetic groups have strong analytical properties arising from the negative holomorphic sectional curvatures in $TD$.  Similar results hold in the general case as a consequence of the horizontality of a VHS.  There are also analogous properties that  hold for the curvatures of the Hodge bundles.   These results are classical; cf.~\cite{C-MS-P} for an exposition and references.
  In this partly expository paper we will not discuss these curvature properties, but rather we will  focus on variations of Hodge structure from an EDS perspective.

  We begin with the
  \begin{demo}\lab{3.1}
  {\bf Basic observation.} {\it Integral elements of $I$ are given by abelian subalgebras}
  \[
  \scrA\subset \scrG^{-1,1} .
  \]
  \end{demo}

This result is a consequence of the Maurer--Cartan equation; it is discussed in a more general context at the beginning of Section~\ref{4} below (cf.\ the proof of Proposition~\ref{4.2}.  Here we are using the notations in Section~\ref{2a} above.  At a point $F\in D$, the Lie algebra $\scrG$ has a Hodge structure constructed from $F$, and \pref{3.1} above identif\/ies the space of integral elements of $I$ in $T_F D$.

  In the theory of EDS, there are two types of exterior dif\/ferential systems, which may be informally described as follows:
  \begin{itemize}\itemsep=0pt
  \item {\it Those that may be integrated by ODE methods.}
  \item {\it Those that require PDE techniques}.
  \end{itemize}
  Although the f\/irst type may be said to be ``elementary",  it includes many important EDS's that arise in practice.  Roughly speaking, these are the exterior dif\/ferential  systems that, by using ODE's, may be  put in a standard local normal form.

  We begin by f\/irst recalling the contact system and then giving two interesting elementary examples   that arise in variations of Hodge structures.

  \begin{cons} Here, $\dim M=2k+1$ and $I$ is locally generated by a $1$-form $\theta$ with $\theta\wedge (d\theta)^k\ne 0$.  Then, by using ODE methods, it may be shown that there are local coordinates $(x^1 ,\dots, x^k,y,y_1,\dots,y_k)$ such that $\theta$ may be taken to be given by (using summation convention)
  \[
  \theta= dy-y_i dx^i .
  \]
  Local integral manifolds are of dimension $\leqq k$, and those of dimension $k$ on which $dx^1\wedge \dots \wedge dx^k\ne 0$ are graphs
  \[
(  x^1,\dots,x^k) \to (x^1,\dots , x^k, y(x),\part_{x^1} y(x),\dots,\part_{x^k}y(x)) .\]
  \end{cons}

  \begin{propn} For weight $n=2$ and Hodge numbers $h^{2,0}=2$, $h^{1,1}=k$, $I$ is a contact sys\-tem~{\rm \cite{CT2}}.
  \end{propn}

  \begin{proof} We shall use the proof to illustrate and mutually relate the two main computational methods
  that have been used to study the dif\/ferential system $I$.

  The f\/irst is using the Hodge structure on $\scrG$ to identify both the tangent space at $F\in D$ and the f\/ibre $W_{I,F}$ with subspaces of $\scrG_\C$.

  For weight $n=2$ we have for $\pp \cong T_FD$ the inclusion
  \begin{equation}\lab{3.2}
  \pp \subset \Hom(H^{2,0},H^{1,1})\oplus \Hom(H^{2,0},H^{0,2})\oplus \Hom(H^{1,1},H^{0,2}) .
  \end{equation}
  Using the polarizing form $Q$ gives
  \[
  \left\{
  \begin{array}{l}
  H^{1,1}\cong \checkh^{1,1},\\
  H^{0,2}\cong \checkh^{2,0} .\end{array}\right.\]
  Then $\pp$ is given by $A=A_1\oplus A_2\oplus A_3$ in the RHS of \pref{3.2} where
  \begin{equation}\lab{3.3}
  \left\{
  \begin{array}{l}
  A_2=-{}^t\! A_2,\\
  A_3= {}^t \!A_1 .\end{array}\right.\end{equation}
  The transpose notation refers to the identif\/ication just preceding~\pref{3.2}.
  The sub-space $\scrG^{-1,1}$ is def\/ined by $A_2=0$.  By~\pref{3.1}, integral elements are then given by linear subspaces
  \begin{equation}\lab{3.4}
  \scrA \subset \Hom(H^{2,0},H^{1,1})
\end{equation}
satisfying
\begin{equation}\lab{3.5}
A\, {}^t\!B = B\,{}^t\! A
  \end{equation}
  for $A,B\in \scrA$.

  The second method is via moving frames.  Over an open set in $D$ we choose a holomorphic frame  f\/ield
  \[
  \underbrace{ \underbrace{e_1,\dots e_h}_{F2},f _1,\dots,f_k}_{F^1} , e^\ast_1,\dots,e^\ast_h\]
  adapted to the Hodge f\/iltrations and to $Q$ in the sense that
  \[
  \left\{\begin{array}{l}
  Q(e_\alpha,e^\ast_\alpha)  = 1,\\[4pt]
  Q(f_i,f_i)   = 1\end{array}\right.
  \]
  and all other inner products are zero.  Then, using summation convention,
  \[
  \left\{
  \begin{array}{l}
  de_\alpha    =\theta^\beta_\alpha e_\beta+\theta^i_\alpha f_i+\theta_{\alpha\beta}e^\ast_\beta,\\[4pt]
  df_i   = \theta^\alpha_i e_\alpha + \theta_{ij}f_j+\theta_{i\alpha}e^\ast_\alpha .\end{array}\right.\]
  Denoting $T_FD$ by just $T$ and referring to \pref{3.2}, \pref{3.3} we have
  \begin{gather*}
  \theta^i_\alpha\in T^\ast \otimes \Hom(H^{2,0},H^{1,1}) \leftrightarrow  A_1,\\
  \theta_{\alpha\beta}\in T^\ast\otimes \Hom(H^{2,0},H^{0,2}) \leftrightarrow  A_2 .
  \end{gather*}
  From $0=dQ(e_\alpha,e_\beta)=Q(de_\alpha,e_\beta)+Q(e_\alpha,de_\beta)$ we have
  \[
  \theta_{\alpha\beta}+\theta_{\beta\alpha}=0\]
  which is the f\/irst equation in~\pref{3.3}.  The second equation there is
  \[
  \theta^\alpha_i = \theta^i_\alpha .
  \]
  From the above we see that $I$ is generated by the $h(h-1)/2$ $1$-forms $\theta_{\alpha\beta}$ for $\alpha<\beta$, where $h=h^{2,0}$.  When $h=1$, we are in the classical case and $I$ is zero.  When $h=2$, $I$ is generated by a single $1$-form $\theta=\theta_{12}$ and is thus a candidate to be a contact system.

  To calculate $d\theta$ we use $d(de_1)=0$, which together with the above formulas give
  \[
  d\theta \equiv \theta^i_1\wedge \theta^i_2 \quad \mod \theta .
  \]
  Since the $1$-forms $\theta^i_\alpha$ are independent, we see that $\theta\wedge (d\theta)^k\ne 0$ as desired.\end{proof}

  The simplest case of this example is when $k=1$.  Then $\dim D=3$ and $I$ is locally equivalent to the standard contact system
 \begin{equation}
 \lab{3.7}
 dy-y' dx = 0\end{equation}
 in $\C^3$ with coordinates $(x,y,y')$.

 In general, as noted above the contact system is locally equivalent to the standard one generated by
 \[
 \theta= dy-y_i dx^i \]
 in $\C^{2k+1}$ with coordinates $(x^1,\dots,x^k,y,y_1,\dots,y_k)$.  Their integral manifolds are graphs
 \[
 \big(x^1,\dots,x^k,y(x),\part_{x^1}y(x),\dots,\part_{x^k}y(x)\big)\]
 where $y(x)$ is an arbitrary function of $x^1,\dots,x^k$.  It is of interest to see explicitly how one may construct integral manifolds depending on one arbitrary function of $k$ variables.

 For this it is convenient to choose a basis for $H_\C$ relative to which
 \begin{gather*}
  Q =  \begin{pmatrix}
\enspace  0&0&I\enspace  \\[4pt]
\enspace   0&-I&0 \enspace \\[4pt]
 \enspace  I&0&0\enspace  \end{pmatrix}\hskip-7pt
 \begin{array}{l}
 \} {\scriptstyle h}\\[4pt]
\} {\scriptstyle k}\\[4pt]
 \} {\scriptstyle h}\end{array}\\[-9pt]
\hspace*{12.5mm} {\mbox{\tiny$\underbrace{\phantom{I}}_{h} \ \ \,   \underbrace{}_k \  \  \underbrace{}_h$}}
 \end{gather*}
 Then $F^2$ will be spanned by the columns in a matrix
 \begin{gather*}
 F =  \begin{pmatrix}
\enspace A\enspace\\[4pt]
\enspace B\enspace\\[4pt]
 \enspace C\enspace\end{pmatrix} \hskip-7pt
 \begin{array}{l}
 \} {\scriptstyle h}\\[4pt]
\} {\scriptstyle k}\\[4pt]
 \} {\scriptstyle h}\end{array}\\[-9pt]
\hspace*{11mm} \enspace {\mbox{\tiny$\underbrace{\phantom{I}}_{h}$}}
  \end{gather*}
  In an open set we will have $\det C\ne 0$, and it is convenient to normalize to have $C=I$.  The equations ${}^t \! FQF=0$ are
  \begin{equation}\lab{3.8}
  A+{}^t\!A ={}^t\!BB .
  \end{equation}
  The inf\/initesimal period relation ${}^t\!FQ\, d F=0$
  is
  \begin{equation}\lab{3.9}
  dA={}^t\!B\, dB .\end{equation}
  Note that adding to this relation its transpose gives the dif\/ferential of~\pref{3.8}.  Thus
  \begin{demo}\lab{3.10}
  {\it If $\pref{3.9}$ is satisfied and $\pref{3.8}$ is satisfied at one point, then it is satisfied identically.}
  \end{demo}

  We now specialize to the case $h=2$ and write
  \[
  A= \begin{pmatrix}
  A_{11}&A_{12}\\[4pt]
  A_{21}&A_{22}\end{pmatrix},
  \qquad B= \begin{pmatrix}
  B_1&x^1\\[2pt]
  \vdots &\vdots\\[4pt]
  B_k&x^k\end{pmatrix}, \]
 where  the $A_{\alpha\beta}$, $B_\alpha$, $x^j$ are independent variables and eventually
 the $A_{\alpha\beta}$ and $B_j$ are to be functions of $x^1,\dots,x^k$.  The equations \pref{3.9} are
 \begin{alignat*}{3}
& (i)  &&\qquad   dA_{11}  =  \Sigma B_jdB_j=d\left( \frac{\sum B^2_j}{2}\right) ;& \\
& (ii) &&\qquad dA _{22}  =  \Sigma x^j dx^j = d\lrp{\frac{\sum (x^j)^2}{2}}; & \\
& (iii) &&\qquad dA_{12}  =  B_j dx^j; & \\
& (iv) &&\qquad dA_{21}  =  x^j dB_j.&
 \end{alignat*}
Now $0 = d(dA_{12})=\Sigma d B_j\wedge dx^2$
 or{\samepage
 \[
 B_j(x) = \partial_{x^j} B(x),
 \]
 where $B(x)$ is an arbitrary function of $x^1,\dots ,x^k$.}

 This then leads to the construction of integral manifolds, here constructed by ODE methods.  As previously noted, in general, PDE methods~-- the Cartan--K\"ahler theorem~-- are required.

 In the above example the $1^{\rm st}$ derived system $I_{[1]}=(0)$.  Below we shall give an  example,  usually referred to as the {\it mirror quintic}, that shows a dif\/ferent phenomenon.
  In general we have
    \begin{gather}\lab{3.6}
  d\theta_{\alpha\beta}\equiv \theta^i_\alpha\wedge \theta^i_\beta \quad \mod \mathrm{span}\{\theta_{\alpha\beta}\hbox{'s}\} .\end{gather}
  We may then think of $I$ as a sort of {\it multi-contact system} (cf.\ \cite{M}).  Note that for each $\alpha<\beta$
  \[
  \theta_{\alpha\beta} \wedge (d\theta_{\alpha\beta})^k\ne 0\]
  and
  \[
  \bigwedge_{\alpha<\beta}\lrp{\theta_{\alpha\beta}\wedge (d\theta_{\alpha\beta})^k} \ne 0 .
  \]
  Integral elements are spanned by matrices $A^i_{\alpha\la}$ satisfying (using summation convention)
  \[
  A^i_{\alpha\la}A^i_{\beta\mu}=A^i_{\alpha\mu} A^i_{\beta\la}\]
  corresponding to subspaces $\mathrm{span}\{A^i_{\alpha\la}\}\subset \Hom(H^{2,0},H^{1,1})$ on which the RHS of \pref{3.6} vanishes.

  When $h= h^{2,0}=2$, the maximal abelian subalgebras \pref{3.4} have dimension $k=h^{1,1}$ and correspond to Lagrangian subspaces in the symplectic vector space $\scrG^{-1,1}$.  In general, the multi-contact nature of $I$  suggests a bound of the sort
  \[
  \dim \scrA \leqq \frac{1}{2} (h^{2,0} h^{1,1}) .
  \]
  This is in fact proved in \cite{C}, where it is also shown that the bound is sharp for $h^{1,1}$ even.  When $h^{1,1}$ is odd, it is shown there that the sharp bound is $\frac12 h^{2,0}(h^{1,1}-1)+1$.
 Because it illustrates yet another method for estimating the dimension of, and actually constructing, integral elements we shall give a proof of this result in the cases $h^{2,0}=2$ and $h^{2,0}=3$.  By \pref{3.5}, we are looking for a linear subspace
 \[
 E\subset \Hom(H^{2,0},H^{1,1})\]
 with a  basis $A_1,\dots,A_r$ satisfying
 \[
 A_i {}^t\!A_j=A_j{}^t\!A_i\]
 for all $i,j$.  A key fact is
 \[
 A_i(v)=0\Rightarrow Q (A_j(v), \hbox{Image}\, A_i)=0\]
 for all $j$.

 \subsection*{Case $\boldsymbol{h^{2,0}=2}$}

 We want to show that $r\geqq k+1$ is impossible.  Assuming that $r\geqq k+1$ then some $A\in E$ drops rank.  With suitable labelling, for some $0\ne v\in H^{2,0}$ we will have $A_{k+1}(v)=0$.  By the above, this gives
 \[
 A_i(v)\in \rim(A_{k+1})^\bot \]
 for $1\leqq i\leqq k$.
 Remembering that $h^{1,1}=k$, we have  $\dim\rim(A_{k+1})^\bot\leqq k-1$ so that $A_1(v),\dots,A_k(v)$ are linearly dependent.  Relabelling, we may assume that $A_k(v)=0$.  Now
 \[
 \dim(\rim A_k+\rim A_{k+1})\geqq 2 ,\]
 since otherwise some linear combination of $A_k$ and $A_{k+1}$ is zero.  Then
 \[
 A_i(v)=(\rim A_k+\rim A_{k+1})^\bot\]
 for $i=1,\dots,k-1$ forces a linear dependence on $A_1(v),\dots,A_{k-1}(v)$.  Proceeding by downward induction gives a contradiction to the assumption $r\geqq k+1$.

 \subsection*{Case $\boldsymbol{h^{2,0}=3}$}

 Since $r\geqq k-1$ we have $A(v)=0$ for some $0\ne v\in H^{2,0}$, $A\in E$.  Choosing a basis $A_1,\dots,A_m$ for the kernel of the map $v\to A(v)$, we let
 \[
 s=\dim \{\rim A_1+\dots+\rim A_m\} .\]
 This gives an injective map
 \[
 H^{2,0}/ \C v\to \C^s\subset \C^k .\]
 By the previous case, $m\leqq k$.  Clearly $m\leqq 2s$ since $\dim(\Hom H^{2,0}/ \C v,\C^k)=2s$.  Thus $m\leqq \min (2s,k)$.  Since for all $A\in E$
 \[
 A(v)\in (\rim A_1+\dots+\rim A_m)^\bot \cong\C^{k-s}\]
 we have $m\geqq r-(k-s)$.  Thus $\min (2s,k)\geqq r-(k-s)$ which gives
 \[
 k-s+\min(2s,k)\geqq r\]
or $r\leqq 3k/2$ as desired.

 When $k=2s$ this bound is sharp.  To see this, take a subspace $\scrU\subset H^{1,1}\cong \C^{2s}$ such that $Q|_\scrU=0$.  Let $e_1,e_2,e_3\in H^{2,0}$ be a basis, and take linearly independent $A_1,\dots,A_k$ satisfying
 \[
 e_1\to 0\qquad \hensp{and} \qquad e_2,e_3\to\scrU .\]
 Then take linearly independent $B_1,\dots,B_s$ satisfying
 \[
 e_1\to\scrU \qquad \hensp{and} \qquad e_2,e_3\to 0 .\]
 It follows that  $E=\mathrm{span}\{A_1,\dots,A_k,B_1,\dots,B_s\}$ is an integral element of dimension $k+s=3k/2$.

 \begin{exam} We take weight $n=3$ and Hodge numbers $h^{3,0}=h^{2,1}=1$.  We will show that in this case $\dim D=4$ and $I$ is locally equivalent to the Pfaf\/f\/ian system
 \begin{equation}\lab{3.11}
 \left\{
 \begin{array}{l}
 dy-y'dx  = 0,\\[4pt]
 dy'-y'' dx  = 0 \end{array}\right.\end{equation}
 in $\C^4$ with coordinates $(x,y,y',y'')$.  Since $d(dy-y'dx)= -dy'\wedge dx = -(dy' - y''dx)\wedge dx$,
  the $1^{\rm st}$ derived system of \pref{3.11} has rank one.  In fact, \pref{3.11} is just the $1^{\rm st}$ prolongation of the contact system~\pref{3.7}.
  \end{exam}

  Anticipating the later discussion of the period domain $D$ for general weight $n=3$ Hodge structures with $h^{3,0}=1$ and $h^{2,1}=h$, locally over an open set $\scrU$ we consider a holomorphic frame f\/ield
  \[
  \underbrace{
  	\underbrace{
		\underbrace{e_0}_{F^3}, e_1,\dots , e_h}_{F^2}, e^\ast_1,\dots , e^\ast_h}_{F^1}
, e^\ast_0\]
relative to which{\samepage
\begin{equation}\lab{3.12}
\left\{
\begin{array}{l}
Q(e_0,e^\ast_0)  = 1=-Q(e^\ast_0,e_0),\\[4pt]
Q(e_\alpha,e^\ast_\alpha)   = 1=-Q(e^\ast_\alpha,e_\alpha)\end{array}\right.\end{equation}
and all other pairings are zero.  We observe that}

{\it The system $I$ is equivalent to each of the following}
\[
\left\{
\begin{array}{l}
de_0\equiv 0 \quad  \mod F^2,\\[4pt]
Q(de_0,F^2)=0 ,\end{array}\right.\]
the equivalence resulting from $F^2=(F^{2})^{\bot}$.

 We set, again using summation convention,
 \begin{equation}\lab{3.13}
 \left\{
 \begin{array}{ll}
 de_0   \equiv \theta^\alpha e_\alpha + \theta_\alpha e^\ast_\alpha + \theta e^\ast_0 & \mod F^3,\\[4pt]
 de_\alpha   \equiv \theta_{\alpha\beta}e^\ast_\beta + \theta^\ast_\beta e^\ast_0&  \mod F^2. \end{array}\right.\end{equation}
 Using these equations and the exterior derivatives of \pref{3.12} we have
 \[
 \left\{\begin{array}{l}
 \theta_{\alpha\beta}  = \theta_{\beta\alpha},\\[4pt]
 \theta^\ast_\alpha  = \theta_\alpha .\end{array}\right.\]
 We may then conclude
 \begin{demo}\lab{3.14}
 {\it The $1$-forms $\theta$, $\theta^\alpha$, $\theta_\alpha$, $\theta_{\alpha\beta}$ for $\alpha\leqq \beta$ are semi-basic for the fibering $G_\C\to \checkd$, and over~$\scrU$ give bases for the cotangent spaces.}
 \end{demo}

 \noindent In terms of the Lie algebra description of the tangent space we have
 \[
 \left\{
 \begin{array}{l}
 \theta^\alpha,\theta_{\alpha\beta}   \leftrightarrow \scrG^{-1,1},\\[4pt]
 \theta_\alpha   \leftrightarrow \scrG^{-2,2},\\[4pt]
 \theta  \leftrightarrow \scrG^{-3,3} .\end{array}\right.\]
  From the exterior derivatives of \pref{3.13} we obtain
  \begin{equation}\lab{3.15}
  \left\{
  \begin{array}{l}
  d\theta  = 2\theta^\alpha\wedge \theta_\alpha,\\[4pt]
  d\theta_\beta  = \theta^\alpha\wedge \theta_{\alpha\beta}.\end{array}\right.\end{equation}
  From the Lie algebra description of the tangent space, it follows that $I$ is generated by the Pfaf\/f\/ian equations
  \begin{equation}\lab{3.16}
  \left\{\begin{array}{l}
  \theta  = 0,\\[4pt]
  \theta_\alpha   = 0.\end{array}\right.
  \end{equation}
  From this we see that the $1^{\rm st}$ derived system $I_{[1]}$ is generated by $\theta$, and the $2^{\rm nd}$ derived system $I_{[2]}  = 0$.

  When $h=1$ we are on a $4$-dimensional manifold with a rank 2 Pfaf\/f\/ian system $I$ and where rank $I_{[1]}=1$, $I_{[2]}=0$.  It is well known (Engel's theorem) and elementary~\cite{IL} and  \cite[Chapter~2]{BCGGG}
  that such a system is locally equivalent to \pref{3.11}.

  We now return to the general case when $h^{3,0}=1$ and $h^{2,1}=h$.  Because  of their importance in algebraic geometry, we shall be interested in integral manifolds $S\subset D$ where the map
  \[
  T  S \to \Hom(H^{3,0} ,H^{2,1} )\]
  is an isomorphism.  We shall call these {\it integral manifolds of Calabi--Yau type.}

  \begin{demo} \lab{3.17}
  {\bf Proposition} (cf.~\cite{B-G}). {\it The EDS for integral manifolds of Calabi--Yau type is canonically the $1^{\rm st}$ prolongation of a contact system.}
  \end{demo}

Thus, locally these integral manifolds depend on one arbitrary function of $h$ variables.

  \begin{proof} We set $P=\P H_\C\cong \P^{2h+1}$.  Denoting points of $P$ by homogeneous coordinate vectors $[z] = [z_0,\dots,z_{2h+1}]$, we consider on $H_\C$ the $1$-form
  \[\theta = Q(dz,z) .\]
  Rescaling locally    by $z\to fz$ where $f$ is a non-vanishing holomorphic function we have
  \[
  \theta\to f^2 \theta .\]
  It follows that $\theta$ induces a $1$-form, def\/ined up to scaling, on $P$.  Choosing coordinates so that
  \[
  Q= \begin{pmatrix}
  &&&0& -1\\[4pt]
  &&& 1&0\\[4pt]
  &&\cdot\enspace\raise6pt\hbox{$\cdot$}\enspace\raise12pt\hbox{$\cdot$}&&\\[4pt]
  0&-1&\\[4pt]
  1&0\end{pmatrix}  , \]
  in the standard coordinate system on $P$ where $z_0=1$ we have{\samepage
  \begin{gather*}
  \theta =  dz_1 + \sum^m_{j=1} z_{2j}dz_{2j+1}-z_{2j+1}dz_{2j},\\
  d\theta =  2\sum^m_{j=1} dz_{2j}\wedge dz_{2j+1},
  \end{gather*}
  from which it follows that $\theta$ induces a contact structure on $P$.}

  Now let $S\subset D$ be an integral manifold of Calabi--Yau type and consider the diagram
  \[
  \xymatrix{
  S \hskip-46pt&{\subset}\hskip-26pt&D\ar[d]^{\pi}\\
  & &P}
  \]
  where $\pi(F)=F^3$ viewed as a line in $H_\C$.  Choose local coordinates $s_1,\dots , s_m$ on $S$ and write~$\pi|_S$~as
  \[
  s\to [z(s)] .\]
  Letting $F^p_s\subset H_\C$ denote the subspace corresponding to $s\in S$,  the condition $dF^3_s \subset F^2_s=(F^{2}_s)^{\bot}$ gives
  \[
  \pi^\ast \theta\mid_S = Q(dz(s),z(s))=0  ; \]
  i.e., $\pi(S)$ is an integral manifold of the canonical contact system on $P$.  In terms of $z(s)$, the integral manifold $S$ is given by
  \begin{equation} \lab{3.18}
  \left\{ \begin{array}{l}
  F^3_s   = [z(s)],\\[4pt]
  F^2_s   = {\rm span}\{z(s),\part_{s_1} z(s),\dots,\part_{s_n}z(s)\}\in F(1,h),\end{array}\right.\end{equation}
  where $F(1,h)$ is the manifold of f\/lags $F^3\subset F^2\subset  H_\C$ where $\dim F^3=1$, $\dim F^2=h+1$, and where we have used the obvious inclusion $\checkd \subset F(1,h)$.

  This process may locally be reversed.  Given an $h$-dimensional integral manifold $s\to [z(s)]$ of the canonical Pfaf\/f\/ian system on $P$, we may def\/ine an integral manifold of $I$ by \pref{3.18}.

  Finally, referring to \pref{3.16} we see that what we denoted by $\theta$ there corresponds to the $\theta$ just above.  From equation~\pref{3.15} we may infer that $I$ is in fact locally just the $1^{\rm st}$ prolongation of the canonical contact system on $P$.
\end{proof}

 \begin{rem} In terms of an arbitrary function $g(s_1,\dots s_h)$ of $h$ variables, in the standard coordinate system the above integral manifold of the contact system is
\[
  (s_1,\dots , s_n)\to (g(s), \part_{s_1}g(s),s_1,\dots \part_{s_n} g(s),s_n).
\]
  It follows that  $S\subset D$ is given parametrically in terms of $g(s)$ and its f\/irst and second derivatives.
  \end{rem}

  \subsection[A Cartan--K\"ahler example and a brief guide to some of the literature]{A Cartan--K\"ahler example and a brief guide to some of the literature}
  \lab{3.b}

  \begin{exam}
  We shall now discuss how the Cartan--K\"ahler theory applies to the case of weight $n=2$, $h^{2,0}=3$ and $h^{1,1}=h$.  This is the f\/irst ``non-elementary" case.
  \end{exam}

  The following  is a  purely linear algebra discussion dealing with the data
  \begin{itemize}\itemsep=0pt
  \item Complex vector spaces $P$, $R$ of dimension $3$, $h$ respectively and where $R$ has a non-degenerate symmetric form giving an identif\/ication $R\cong \check{R}$.
  \item In $\Hom (P,R)\cong \check{P}\otimes R$ we are looking for an {\it abelian subspace} $E$, that is  a linear subspace of $\check{P}\otimes R$ such that
  \[
  {}^t\! A B-{}^t\!BA = 0\]
  for all $A,B\in E$.  Here, the transpose is relative to the identif\/ication $R\cong \check{R}$.
  \end{itemize}
  We seek $h$-dimensional abelian subspaces $E$ that are in general
  position relative to the tensor-product structure in the sense that for a general $u\in P$ the composite map
  \[ \begin{array}{ccc}
  E\to &\Hom(u,R)&\to R\\
  &\cap&\\
  &\Hom(P,R)\end{array}\]
  is an isomorphism.
  With this assumption we shall be able to establish that general integral elements that satisfy it are ordinary in the EDS sense.

  It is convenient to choose a basis $e_i$ for $P$ and an orthonormal basis $e_\alpha$ for $R$ so that
  \[
  e_{i\alpha}=: \check{e}_\alpha \otimes e_i
  \]
  gives a basis for $\Hom(P,R)$.  Denote by
  \[
  \theta^\alpha_i\in\Hom(P,R)^{\vee}\]
  the dual basis.  The notation has  been chosen to line up with the ``moving frame" notations above.  Using the summation convention and setting
  \[
  \Om_{ij}= \theta^\alpha_i\wedge \theta^\alpha_j=-\Om_{ji},
  \]
  abelian subspaces are given by linear subspaces $E$  on which the restriction
  \begin{equation}\lab{3.20}
  \Om_{ij}\mid_E=0.
  \end{equation}

  We may assume without loss of generality that $E$ is of dimension $h$ and that the condition of general position is satisf\/ied for $u=e_i$, $i=1,2,3$.
  Setting $\Theta_i=\wedge_\alpha \theta^\alpha_i$, we have
    \begin{equation}\lab{3.21}
  \Theta_i\mid_E\ne 0.
  \end{equation}
  Then $E$ is  def\/ined by linear equations
  \[
  \left\{
  	\begin{array}{ll}
	\theta^\alpha_2   =A^\alpha_\beta\theta^\beta_1,&\qquad \det \| A^\alpha_\beta\|\ne 0,\\[4pt]
	\theta^\alpha_3   = B^\alpha_\beta\theta^\beta_1,&\qquad \det \|B^\alpha_\beta\|\ne 0
	 .\end{array}\right.\]
The equations \pref{3.20} for $\Om_{12}$ and $\Om_{13}$ are (Cartan's lemma)
\[
\left\{
\begin{array}{l}
A= {}^t\!A,\\[5pt]
B={}^t\!B .\end{array}\right.\]
The equation \pref{3.20} for $\Om_{23}$ is
\begin{equation}\lab{3.22}
[A,B^{-1}]=0 .\end{equation}
In fact, if $\theta^\alpha_1=C^\alpha_\beta\theta^\beta_3$ on $E$, then $A=CB$ and $C={}^t C$ gives \pref{3.22}.

\begin{demo}\lab{3.23}
{\bf Proposition.} {\it The dimension of the space of generic $h$-dimensional abelian subspaces is $h(h+3)/2$.}
\end{demo}

\begin{proof} The group operating on $\Hom(P,R)$ is $\GL(P)\times O(R)$; this group preserves the set of equations \pref{3.20} that def\/ine the abelian subspaces of $\Hom(P,R)$.  For elements $I\times T$ in this group, the action on the matrices $A$ and $B$ above is given by
\[
\left\{ \begin{array}{l}
	A   \to {}^t  TAT,\\
	B   \to \upt TBT .\end{array}\right.\]
 For $A$ generic we may f\/ind $T$ so that $A={\rm diag}(\la_1,\dots,\la_h)$ where the $\la_\alpha$ are non-zero and distinct.  The condition \pref{3.22} then gives that $B={\rm diag} (\mu_1,\dots,\mu_h)$.  We thus have as a basis for~$E$
 \[
 v_\alpha=e_{1\alpha}+\la_\alpha e_{2\alpha}+\mu_\alpha e_{3\alpha}, \]
 where $1\leqq \alpha\leqq h$.
 \end{proof}

 \begin{demo}\lab{3.24}
 {\bf Proposition.} {\it The equations that define the polar space $H(E)$ have rank $2h$.  Thus $E$ is a maximal integral element.}
 \end{demo}

\begin{proof}
 For $v=C^\alpha e_{1\alpha}+D^\alpha e_{2\alpha}+E^\alpha e_{3\alpha}$ the equations that express the condition that $\mathrm{span} \{E,v\}$ be an integral element are
 \[
 \Om_{ij}(v,v_\beta)=0 .\]
 These compute out to be, {\it with no summation on repeated indices},
 \[
 \left\{\begin{array}{l}
 	D^\beta   = \la_\beta C^\beta,\\[4pt]
	E^\beta  = \mu_\beta C^\beta,\\[4pt]
	\mu_\beta D^\beta   = \la_\beta E^\beta .\end{array}\right.\]
These are compatible and have rank $2h$.  For any solution, we replace $v$ by $\tilde{v}= v-C^\alpha v_\alpha$ (summation here), and then $\tilde v=0$ by the above equations when all $C^\alpha = 0$.
\end{proof}

\begin{demo}\lab{3.25}
{\bf Proposition.}  {\it A general abelian subspace $E$ that is in general position with respect to the tensor product structure is an ordinary integral element.}\end{demo}

\begin{proof} In order to not have the notation obscure the basic idea, we shall illustrate the argument in the case $h=2$.  A general line $E_1$ in $E$ is spanned by $w=\rho^\alpha v_\alpha$ (summation convention).  The polar equations for $v$ as above are $\Om_{ij}(w,v)=0$ for $i<j$.  These are
\[
\left\{
\begin{array}{l}
	\rho^1 (D^1- \la_1 C^1)+ \rho^2(D^2-\la_2 C^2)     = 0,\\[5pt]
	\rho^1(E^1-\mu_1 C^1)+\rho^2(E^2-\mu_2 C ^2)  =0, \\[5pt]
	\rho^1(\la_1 E^1-\mu_1D^1)+\rho^2(\la_2 E^2-\mu_2 D^2)  =0 .\end{array}\right.
\]
For a general choice of $\rho^\alpha$, $\la_\alpha$, $\mu_\alpha$ these are independent.  Hence, in the notation of Section~\ref{2.b}
we have
\[
c_0=0,\qquad c_1=3 .\]

On the other hand, the $2$-dimensional abelian subspaces in $\Hom(P,R)\cong \C^6$ is a smooth, $5$-dimensional subvariety in the $8$-dimensional Grassmannian $\Gr(2,6)$.  Therefore the codimension of $E$ in $\Gr(2,6)$ is $3=8-5$.  Consequently, Cartan's test is satisf\/ied.
\end{proof}

In \cite{CT} it is shown how to integrate the above system, with the result
\begin{demo}\lab{3.26}
{\it Integral manifolds of the above system are parametrized by generating functions $f_1$, $f_2$ subject to the PDE system}
\[
[H_{f_1},H_{f_2}]=0,
\]
where $H_f$ is the Hessian matrix of $f$.
\end{demo}

\subsection*{A brief guide to some of the literature}

Building on earlier works, in \cite{M} a number of results are proved:
\begin{itemize}\itemsep=0pt
\item[$(i)$]  Under rather general assumptions of the type
\begin{itemize}\itemsep=0pt
\item[$\bullet$] the weight $n=2m+1$ is odd, $h^{m+1,m}>2$ and all other $h^{p,q}>1$,
\item[$\bullet$] the weight $n=2m$ is even, $h^{m,m}>4$, $h^{m+2,m-2}>2$ and all other $h^{p,q}>1$
\end{itemize}
a bound on the dimension of integral elements is obtained, and it is shown that there is a~unique integral element $E$ that attains this bound.
\item[$(ii)$] There is a~unique germ of integral manifold whose tangent space is $E$.
\end{itemize}
The following interesting question arises:
\begin{quote}
{\it Is $E$ an ordinary integral element?}
\end{quote}
This is a question of computing the ranks of the polar equations of a generic
f\/lag in $E$.  Setting $\dim E=p$ and $\dim D=N$, since $E$ is unique we have at $E$
\[
\codim G_p(I)=\dim \Grass(p,N)=p(N-p)\]
and Cartan's test is whether the inequality
\[
p(N-p)\geqq c_1+\dots+c_{p-1}\]
is an equality.  Suppose that equality does hold, so that $E$ is ordinary and the Cartan--K\"ahler theorem applies.  As explained in \cite{IL} and \cite{BCGGG} there are {\it Cartan characters} $s_0,s_1,\dots,s_p$ expressed in terms of the $c_i$ such that local integral manifolds of~$I$ depend on $s_p$ arbitrary functions of $p$-variables, $s_{p-1}$ arbitrary functions of $(p\!-\!1)$-variables, $\dots,s_0$ arbitrary functions of $0$-variables (i.e., constants).  Mayer's result would then follow by showing that equality holds in Cartan's test and
\[
s_p=
\dots=s_1=0,\qquad s_0=1 .\]
In fact, if $E$ is ordinary then the Cartan characters must be given in this way.
Another interesting question is
\begin{quote}
{\it Do the Mayer integral manifolds arise from algebraic geometry?}
\end{quote}
Since there is a unique  one of these through each point of $D$, one knows that in the non-classical case a general Mayer integral manifold does not arise from geometry.  In \cite{CS} it is shown that in weight two with $h^{1,1}$ even {\it some} maximal integral manifolds are realized geometrically.

Among the works that preceded \cite{M}, and in some cases led up to it are:
\begin{enumerate}\itemsep=0pt

\item[\cite{CT2}:] In this paper it was recognized that the weight two horizontal distribution is a generalization of the contact distribution.

\item[\cite{CD}:] Here it is shown that most hypersurface variations are maximal, i.e.\ their tangent space at a point is not contained in a larger integral element.  These integral elements are not of maximal dimension in the sense of~\cite{M}.

\item[\cite{CKT}:] This paper gives the results in weight two that were extended to the general case in~\cite{M}.
\end{enumerate}

\smallbreak We conclude this section by   discussing the recent work \cite{A}.

  One may think of an integral element $E$ at $F\in D$ as giving an action of $\Sym E$ on  $\opplus_p F^p/F^{p+1}$ that is compatible with $Q$; i.e., the action preserves the pairings $Q:F^p/F^{p+1}
\!\otimes\! F^{n-p}/F^{n-p+1}{\to} \C$.  Especially if $\dim E$ is not small, such an action will have many algebraic invariants.  Experience suggests that those that  arise from geometry will have a special structure.

One such relates to the notion of {\it symmetrizers} due to Donagi (cf.~\cite{D}
and~\cite{CD}) and used in~\cite{A}.  Given vector spaces $A$, $B$, $C$ and a bilinear map
\[
\Phi: \ A\otimes B\to C\]
one def\/ines
\[
\Sym\Phi=\lrc{
	\Psi\in \Hom(A,B):\Phi(a,\Psi(a'))=\Phi(a',\Psi(a))\hensp{for}a,a'\in A}.	\]
One may consider the $\Sym$'s of the various maps
\[
\Sym^kE\otimes H^{p,n-p}\to H^{p-k,n-p+k}
\]
arising above.  In \cite{A} these are studied when $\dim E$ and the $h^{p,n-p}$ are the same as for smooth hypersurfaces $X\subset \P^{n+1}$ where $n\geqq 3$ and $\deg X\geqq n+3$, and where the special structure arising from Macauley's theorem is satisf\/ied.  It is shown that a particular Sym is generically zero, but is non-zero in the geometric case.  Thus the integral elements arising from hypersurface deformations, which are known to be maximal, satisfy non-trivial algebraic conditions and are thus non-generic.

 We shall give a brief discussion of proofs of these results, based on the paper \cite{Gr}, in the case $n=2$.  For this we use the notations from the proof of Theorem~\ref{5.10} below.  For smooth surfaces in $\P^3$ the analogous identif\/ication to \pref{5.13} is
 \[
 H^{2-p,p}\cong V^{(p+1)d-4}/J_{(p+1)d-4} .\]
 We also have an identif\/ication
 \[
 E\cong V^d / J_d .\]
 Using $E\subseteq \Hom(H^{2,0},H^{1,1})$ there is a map
 \[
 \Hom (H^{2,0},E)\to \Hom(\ottimes^2 H^{2,0},H^{1,1}) .\]
 From the above identif\/ications, we have
 \[
 V^4/J_4\subseteq \Hom(H^{2,0},E)\cong \Hom(V^{d-4},V^{2d-4}/J_{2d-4})\]
 and under the above map $V^4/J_4$ lands in $\Hom(\Sym^2 H^{2,0},H^{1,1})$.  Put another way
 \[
 V^4/J_4\subseteq \ker\lrc{
 	\Hom(H^{2,0},E)\to \Hom(\wedge^2 H^{2,0},E)}.
\]
It may be seen that for $d\geqq 5$ and for $E$ a general integral element, this kernel is non-zero.
This is illustrative of the non-genericity results in~\cite{A}.

If follows from \cite{Gr} that for $d\geqq 5$
\[
V^4 \to \check{H}^{2,0}\otimes E\to \wedge^2 \check{H}^{2,0}\otimes H^{1,1}\]
is exact at the middle term.  Thus one has

\medbreak\noindent {\bf Integrability.} The composition
\[
H^{2,0}\to \check E\otimes H^{1,1}\to \wedge^2 \check E\otimes H^{0,2}\]
is zero.

\medbreak\noindent {\bf Non-generic condition.} The kernel
\[
\ker\lrc{
	\check H^{2,0}\otimes E\to \wedge^2 \check H^{2,0}\otimes H^{1,1}
	}
\]
is non-zero, and as a vector space is isomorphic to $V^4$.

One may next proceed, for $d\geqq 8$, to look at
\[
V^{d-8} \to \Hom(V^4,H^{2,0})\to \Hom (\wedge^2 V^4,E) ,\]
and by \cite{Gr} this will be exact at the middle term of $d\geqq 9$.
This process continues and leads to a large number of non-linear constraints on the integral element arising from variations of a~smooth surface of degree $d\geqq 5$.

Additional constraints, of a dif\/ferent character, arise as follows.  Recall that polynomials $f_0$, $f_1$, $f_2$, $f_3$, where $f_i\in V^d$, form a {\it regular sequence} if $\mathrm{Var}(f_0,f_1,f_2,f_3)=\emptyset$ in $\P^3$.  A special case is when
\[
f_i(x)=F_{x_i}(x)\]
where $F(x)=0$ def\/ines a smooth surface in $\P^3$.  The symmetrizer constructions and results in~\cite{Gr} illustrated above work equally well for a regular sequence.  Thus, for the integral element def\/ined by a surface variation, there is the additional constraint that the ideal generated by a~regular sequence be a Jacobian ideal.

Although we shall not give the argument (using, as usual, Macaulay's theorem), it may be proved that an integral element $E$ def\/ined by a regular sequence is maximal.

We set $\dim E=p$ and denote by
\[
G_p(\hbox{surfaces}) \subset G_p(RS)\subset G_p(I)\]
the integral elements corresponding respectively to smooth surfaces, to regular sequences in $\P^3$, and general integral elements of the dif\/ferential ideal generated by $I$.  As just noted, the f\/irst inclusion is strict.  The following would seem to be an interesting question:
\begin{quote}
{\it Is an integral element $E\in G_p(RS)$ ordinary $($cf.~$\pref{2.11}$ above$)$?}
\end{quote}

\subsection{The derived f\/lag of the EDS associated to a VHS}
\lab{3.c}

We shall prove the

\begin{demo}\lab{3.26n}
{\bf Theorem.}  {\it Let $D$ be a period domain for polarized Hodge structures of weight $n\geqq 2$ and where all $h^{p,q}\ne 0$.  Denote by $I\subset T^\ast D$ the Pfaffian system given by the infinitesimal
period  relations.  Then for the derived system}
\[
I_{[m]} = 0\qquad \hensp{for} \quad m\geqq \frac{\log n}{\log 2} .\]
\end{demo}

\begin{cor} When $n\geqq 2$, $I$ has no completely integrable subsystem.
\end{cor}

\begin{rems}
$(i)$ This corollary roughly means that {\it there are no finite equations satisfied by integral manifolds of $I$.}  More precisely, locally there are no non-constant functions $f$ on open subsets $\scrU$ of $D$ such that, for a point $F\in \scrU$, all local integral manifolds of $I$ passing through $F$ lie in the level set $f=\hbox{const}$.

$(ii)$ The theorem actually is valid on the dual classifying space $\check{D}$.  Over $\check{D}$ there are Pfaf\/f\/ian systems $I(k)$ def\/ined for $k\geqq 1$ by
\begin{equation}\lab{3.27n}
dF^p\subseteq F^{p-k} .\end{equation}
For the VHS Pfaf\/f\/ian system, $I=I(1)$.  We shall establish a stronger version of the theorem by showing that:

{\it If all $h^{p,q}\ne 0$, then}
\begin{equation}
\lab{3.28n}
I_{[m]}\subseteq I(2^m) .\end{equation}
{\it Moreover, if $n\geqq 2$ and all $h^{p,q}\geqq 2$ then equality holds in~$\pref{3.28n}$.}

$(iii)$ We shall also see that the above result is sharp.  The f\/irst case where some $h^{p,q}=0$ but $I\ne T^\ast D$ occurs when $n=4$ and $h^{4,0}\ne  0$, $h^{3,1}\ne 0$ but $h^{2,2}=0$.  Then $I_{[\infty]}\ne 0$.  The f\/inite equations satisf\/ied by a VHS may be described as follows:  Denote by $ \Gr(f^2,H_\C)$ the Grassmanian of $f^2$-dimensional planes in $H_\C$.  Then, since $dF^3\equiv 0 \mod F^3$,  the connected integral manifolds of $I$ lie in the f\/ibres of the obvious map
\[
\check D\to \Gr(f^2,H_\C) .\]

$(iv)$ The simplest interesting case where equality fails to hold in \pref{3.28n} is when $n=4$ and $h^{4,0}=1$.  We shall see below that this is ref\/lected in the structure of the derived f\/lag.  The conclusion will be: {\it Among $4$-folds with $h^{4,0}\ne 0$, Calabi--Yau varieties are distinguished by the behaviour of the derived flag associated to the infinitesimal period relation.}
\end{rems}

\begin{proof}
In order not  to have the notational complexity obscure the basic idea, we shall give the argument in the case $n=4$.  From this we hope that the pattern for the general case will be evident.  The assumption that all $h^{p,q}\ne 0$ is equivalent to the positive dimensionality of all graded quotients in the f\/iltrations
\[
F^4\subset F^3\subset F^2\subset F^1\subset F^0=H_\C\]
corresponding to a point of $\check D$.  We consider adapted frame f\/ields
\[
\underbrace{
	\underbrace{
		\underbrace{
			\underbrace{e_{0\alpha}}_{F^4}
		,e_{1\alpha}}_{F^3}
	, e_{2\alpha}}_{F^2}
, e_{3\alpha}}_{F^1}
,e_{4\alpha}\]
satisfying
\begin{equation}\lab{3.29n}
Q(e_{i\alpha},e_{j\beta})= \left\{\begin{array}{ll}
1&\hensp{if} i+j=4,\ \alpha=\beta,\\[4pt]
0&\hensp{otherwise.}\end{array}\right.\end{equation}
As usual we write, using summation convention,
\[
d e_{i\alpha}=\theta^{j\beta}_{i\alpha}e_{j\beta},\]
where
\begin{equation}\lab{3.30}
d\theta^{j\beta}_{i\alpha}=\theta^{k\gamma}_{i\alpha}\wedge \theta^{j\beta}_{k\gamma} .\end{equation}
The forms $\theta^{j\beta}_{i\alpha}$ for $j\geqq i+1$ are semi-basic for the f\/ibering $G_\C\to \check D$; it is only expressions in these that will have intrinsic meaning in the f\/inal formulas.  We will use the notation $\equiv$ for congruence modulo the $\theta^{j\beta}_{i\alpha}$ for $j\leqq i$.
The system $I$ is given by
\begin{equation}\lab{3.31}
I=\{\theta^{j\beta}_{i\alpha}:j\geqq i+2\},
\end{equation}
where the brackets are used to denote the span.

The exterior derivative of~\pref{3.29n} gives
\begin{equation}\lab{3.32}
\theta^{4-j\beta}_{i\alpha} + \theta^{4-i\alpha}_{j\beta}=0 .\end{equation}
For later reference we write this out in detail
\begin{equation}\lab{3.32n}
\left\{\begin{array}{ll}
\theta^{4\beta}_{0\alpha}+\theta^{4\alpha}_{0\beta}=0,&\qquad \theta^{3\beta}_{1\alpha}+\theta^{3\alpha}_{1\beta}=0,\\[4pt]
\theta^{3\beta}_{0\alpha}+\theta^{4\alpha}_{1\beta}=0,&\qquad \theta^{2\beta}_{1\alpha}+\theta^{3\alpha}_{2\beta}=0,\\[4pt]
\theta^{2\beta}_{0\alpha}+\theta^{4\alpha}_{2\beta}=0,\\[4pt]
\theta^{1\beta}_{0\alpha}+\theta^{4\alpha}_{3\beta}=0 .\end{array}\right.\end{equation}
The $1$-forms $\theta^{j\beta}_{i\alpha}$ are all non-zero and are linearly independent modulo these relations; the f\/irst statement uses the assumption that all $h^{p,q}\ne 0$.
From \pref{3.30}, for $j\geqq i+1$
\begin{equation}
\lab{3.33}
d\theta^{j\beta}_{i\alpha}\equiv \sum_{\left\{
{i<l<j\atop\gamma\hfill}  \right.} \theta^{l\gamma}_{i\alpha}\wedge \theta^{j\beta}_{l\gamma}\end{equation}
which by \pref{3.31} gives the inclusion
\[
\{\theta^{j\beta}_{i\alpha}\}_{j\geqq i+3}\subseteq I_{[1]} .\]
Indeed, if $j\geqq i+3$ then if $i<l<j$ we must have either $l-i\geqq 2$, $j-l\geqq 2$ or both.  Now we use \pref{3.31}.
Further, from \pref{3.33}
\[
d\theta^{i+2\beta}_{i\alpha}\equiv \sum_\gamma \theta^{i+1,\gamma}_{i\alpha}\wedge \theta^{i+2\beta}_{i+1\gamma} .\]
Using \pref{3.32n} we may verify that no linear combination of the $\theta^{i+2\beta}_{i\alpha}$ is in $I_{[1]}$.  Thus we have equality in the above inclusion; i.e.\
\[
\lrc{
	\theta^{3\beta}_{0\alpha},\theta^{4\beta}_{0\alpha},\theta^{4\beta}_{1\alpha}}
=I_{[1]} .\]

We next claim that
\[
I_{[2]}=0 .\]
Indeed, denoting by $\uquiv$ congruence modulo the algebraic ideal generated by $I_{[1]}$, by \pref{3.33} we have
\[
\left\{
\begin{array}{l}
d\theta^{3\beta}_{0\alpha}  \uquiv \theta^{1\gamma}_{0\alpha}\wedge \theta^{3\beta}_{1\gamma} + \theta^{2\gamma}_{0\alpha}\wedge \theta^{3\beta}_{2\gamma},\\[4pt]
d\theta^{4\beta}_{0\alpha}  \uquiv \theta^{2\gamma}_{0\alpha}\wedge \theta^{4\beta}_{2\gamma},\\[4pt]
d\theta^{4\beta}_{1\alpha}   \uquiv \theta^{3\gamma}_{1\alpha} \wedge \theta^{4\beta}_{3\gamma}+\theta^{2\gamma}_{2\alpha}\wedge \theta^{4\beta}_{2\gamma}  ; \end{array}\right.\]
no linear combination of the RHS is $\uquiv 0$.
\end{proof}

Denoting by $\theta^j_i$ the matrix $\| \theta^{j\beta}_{i\alpha}\|$, the general pattern is
\begin{alignat*}{3}
& I= \{\theta^j_i:j\geqq i+1\},\quad && \hbox{i.e.}\quad j-i>1=2^0,&\\
& I_{[1]} \subseteq \{\theta^j_i:j\geqq i+3\},\quad  && \hbox{i.e.} \quad  j-i>2=2^1, & \\
& I_{[2]} \subseteq \{\theta^j_i:j\geqq i+5\},\qquad && \hbox{i.e. } \quad  j-i>4=2^2, &\\
& \quad \vdots&&& \\
& I_{[k]}   \subseteq \{\theta^j_i:j-i>2^k\}. &&&
\end{alignat*}
To show that equality holds, we must analyze in general the relations such as \pref{3.32n} that arise from $Q(F^p,F^{n-p+1})=0$.  The analysis is an extension of that given above when $n=4$, the one dif\/ference being the alteration of signs with the parity of $n$.  We shall not write out the argument here.  Rather we shall look more closely at the $n=4$ case when $h^{4,0}=1$ and when $h^{2,2}=0$.

\medbreak\noindent {\bf Example:} $h^{4,0}=1$.  In this case $\theta^{4\beta}_{0\alpha}=0$ (only $\beta=\alpha$ occurs and we have $\theta^{4\alpha}_{0\alpha}=0$ by the f\/irst equation  in \pref{3.32n}), and thus
\[
I\cap I(4)=0 .\]
Conversely, if this relation holds then $h^{4,0}=1$.

\medbreak\noindent {\bf Example:} $h^{2,2}=0$. In this case
\[
\theta^{2\beta}_{i\alpha}=0,\qquad  \theta^{j\beta}_{2\alpha}=0,\qquad i=0,1\quad \hensp{and}\quad j=3,4 .\]
Now
\[
d\theta^{3\beta}_{1\alpha}\equiv \theta^{k\gamma}_{1\alpha} \wedge \theta^{3\beta}_{k\gamma}\]
and thus
\[
I_{[1]} = \lrc{ \theta^{3\beta}_{0\alpha},\theta^{4\beta}_{0\alpha},\theta^{3\beta}_{1\alpha},\theta^{4\beta}_{1\alpha}} .\]
We claim that
\begin{equation}\lab{3.34}
I_{[2]} = I_{[1]} .\end{equation}
This will imply that $I_{[\infty]} = I_{[1]}$
and will verify the  claim that leaves of the foliation given by $I_{[\infty]}$ are the f\/ibres of $\check D\to \Gr(f^2,H_\C)$.

Continuing to denote by $\uquiv$ equivalence modulo the algebraic ideal generated by $I_{[1]}$, from~\pref{3.33} we have
\begin{gather*}
d\theta^{3\beta}_{0\alpha}   \equiv  \theta^{1\gamma}_{0\alpha}\wedge \theta^{3\beta}_{1\gamma}\uquiv 0,\\
d\theta^{4\beta}_{0\alpha} \equiv  \theta^{1\gamma}_{0\alpha}\wedge \theta^{4\beta}_{1\gamma}+\theta^{3\gamma}_{0\alpha}\wedge \theta^{4\beta}_{3\gamma}\uquiv 0 ,
\end{gather*}
and similarly for $d\theta^{3\beta}_{1\alpha}$ and $d\theta^{4\beta}_{1\alpha}$.

We are not aware of a geometric example where $h^{4,0}\ne 0$, $h^{3,1}\ne 0$ but $h^{2,2}=0$.  There are of course geometric examples of  necessarily rigid
CY 3-folds where $h^{3,0}=1$, $h^{2,1}=0$.

\section[Universal cohomology associated  to a homogeneous Pfaf\/f\/ian system]{Universal cohomology associated\\  to a homogeneous Pfaf\/f\/ian system} \lab{4}

\subsection{EDS aspects of homogeneous Pfaf\/f\/ian systems}\lab{3a}

In this section    all Lie groups and Lie algebras will be real and connected.
  Let $G$ be a Lie group, $V\subset G$ a closed subgroup and
  \[
  M=G/V\]
  the resulting homogeneous space.  We denote by $\scrG$, $\fv$ the Lie algebras respecively of $G,V$.

  \begin{defi} A {\it homogeneous Pfaffian system} is given by a $G$-invariant sub-bundle $I\subset T^\ast M$.\end{defi}
  Equivalently, $I^\bot = : W\subset TM$ is a $G$-invariant distribution.  Since $V$ is connected, this in turn is given by an $\fadv$ invariant subspace
  \[
  \fw\subset \scrG/\fv.\]
  If $\tilde{\fw} \subset  \scrG$ is the inverse image of $\fw$  under the projection of $\scrG\to \scrG/\fv$, $\fadv$ invariance is equivalent to $[\fv,\tilde{\fw}]\subseteq \tilde{\fw}$.

  For any subspace $\scrA\subset \fw$, the condition
  \begin{equation}\lab{4.1}
  [\scrA,\scrA]\subseteq \fv\end{equation}
 is well-def\/ined; the inverse image $\tilde{\scrA}\subset \scrG$ of $\scrA$ should satisfy $[\tilde\scrA,\tilde\scrA]\subseteq \fv$.

 \begin{demo}\lab{4.2}
 {\bf Proposition.} {\it The $G$-invariant integral elements of a homogeneous Pfaffian system are given by the $G$-translates of  subspaces $\scrA\subset \fw$ that satisfy $\pref{4.1}$.}
 \end{demo}

\begin{proof}
  For the projection $G\stackrel{\pi}{\to}M$ we have that $\pi^\ast I$ is a $G$-invariant sub-bundle of $T^\ast G$; it is in fact given by the translates of
 \[
\scri = \fw^\bot \subset \check{\scrG} .\]
 For $\theta\in\scri $, the Maurer--Cartan equation gives, for $\zeta,\eta\in\scrG$,
 \begin{equation}\lab{4.3}
 \lra{d\theta,\zeta\wedge\eta} = -\frac12\lra{\theta,[\zeta,\eta]} .\end{equation}
 For the $G$-invariant subspaces $E\subset TM$ given by the translates of $\scrA\subset \fw$, the condition for such a subspace $E$ to be an integral element is therefore
 \[
 d\theta|_E=0\]
 for all $\theta\in I$.  By \pref{4.3} this is equivalent to \pref{4.1}.\end{proof}

 We now suppose that $M$ is a reductive homogeneous space so that there is an $\fadv$ invariant splitting
 \[
 \scrG=\pp \oplus \fv .\]
 Suppose that $\scrA\subset \pp$ is a subspace that gives an abelian Lie-subalgebra; i.e.,
 \begin{equation}\lab{4.4}
 [\scrA,\scrA]=0 .\end{equation}
 In particular, \pref{4.1} is satisf\/ied.  Let $A=\exp\scrA$ be the connected subgroup of $G$ corresponding to~$\scrA$.

 \begin{demo}\lab{4.5}
 {\bf Proposition} (cf.~Theorem 3.15 in~\cite{M}). {\it The orbit $A\cdot (eV)\subset M$ is an integral manifold of $I$.}
 \end{demo}

This is clear: The tangents to the $A$-orbit are the $A$ translates of $\scrA$, and by \pref{4.3} these are all integral elements.

 The space $G_p(I)$ of $p$-dimensional integral elements is acted on by $G_\R$, and by the proposition every $E\in G_p(I)$ satisfying \pref{4.4} is tangent to at least one integral manifold.  This does {\it not} imply either  that $(i)$ $G_p(I)$ is smooth at $E$, or $(ii)$ that if $G_p(I)$ is smooth at $E$, then $E$ is an ordinary integral element; i.e., $I$ may not be involutive at~$E$.

 \subsection[Characteristic cohomology of the homogeneous Pfaf\/f\/ian system associated to a~VHS]{Characteristic cohomology of the homogeneous Pfaf\/f\/ian system\\ associated to a VHS}\lab{3b}

Associated to a dif\/ferential ideal $\scrI$ on a manifold $M$ are its {\it characteristic cohomology groups} $H^\ast_\scrI(M)$, def\/ined to be the de~Rham cohomology  of the complex $(\Om^\bullet(M)/\scrI,d)$ where $\Om^\bullet(M)$ are the  complex valued  smooth forms on $M$.  One may think of $H^\ast_\scrI(M)$ as the {\it universal cohomology groups induced on integral manifolds of $\scrI$.}  We refer to \cite{BG1, BG2} for a general reference as well as some examples of these groups.

In the case of a homogeneous Pfaf\/f\/ian system one may consider the complex $((\Om^\bullet(M)/\scrI)^G,d)$ of $G$-invariant forms.  By standard arguments one has the
\begin{demo}\lab{4.6}
{\bf Proposition.} $H^\ast((\Om^\bullet(M)/\scrI)^G,d)\cong H^\ast(\scrG,\fv;\fw)$.
\end{demo}

Here, $H^\ast(\scrG,\fv;\fw)$ is the Lie algebra cohomology of $(\scrG,\fv)$ relative to $\fw$; everything is taken over the complex numbers.  In the reductive case, which is the one we shall consider, this is def\/ined as follows: First, we have
\[
\scri =\fw^\bot \subset \check{\pp} .\]
To def\/ine the complex $C^\bullet(\scrG,\fv;\fw;\delta)$ we use the map
\[
\scri\to \scri_2\subset \wedge^2\check{\pp}\]
def\/ined by $\theta\to \theta_2$ where, for $\zeta,\eta\in\pp$
\[
\lra{\theta_2,\zeta\wedge \eta} = \lra{\theta,[\zeta,\eta]} .\]
By the Maurer--Cartan equation, up to scaling these are just the $2$-forms $d\theta$ where $\theta\in I$, that together with $I$ generate the dif\/ferential ideal $\scrI$.  Then we set
\begin{equation}\lab{4.7}
C^q(\scrG,\fv;\fw;\delta) = \left\{
\begin{array}{ll}
0,&\hensp{when} q=0,\\[4pt]
(\check{\pp}/\scri )^\fv,&\hensp{when}q=1,\\[4pt]
(\wedge^q \check{\pp}/\wedge^{q-1}\check{\pp}\wedge\scri+\wedge^{q-2}\check{\pp}\wedge\scri_2)^\fv,&\hensp{for} q\geqq 2 .\end{array}\right.\end{equation}
These are just the values at the identity coset $eV$ of the forms in $(\Om^q(M)/\scrI^q)^G$.  The map $\delta$ is  just the exterior derivative, which when written out has the usual form of the dif\/ferential in Lie algebra cohomology.

In applications, one frequently has a discrete subgroup $\Gamma$ acting properly discontinuously on~$M$, and one is interested in the global invariants of maps
\[
f: \ S\to \Gamma\bsl M.\]
There is then an induced map
\begin{equation}\lab{4.8}
f^\ast: \ H^\ast(\scrG,\fv;\fw)\to H^\ast (S) .\end{equation}
We may think of the image of this map as {\it giving the global topological invariants of maps $f:S\to \Gamma\bsl M$ that may be defined universally for all $\Gamma$'s.}

We now turn to the case of a period domain and use the complexif\/ied Lie algebras.   $\scrI$ will be the dif\/ferential ideal generated by $I\subset T^{1,0}D$ and $\ol{I}\subset T^{0,1}D$.  In the classical case when $n=1$, $\fw=0$ and $V\subset G_\R$ is a maximal compact subgroup, and from what is known it may be inferred that
\begin{demo}\lab{4.9}
{\it $H^\ast(\scrG,\fv;\delta)\cong (\wedge^\ast\check{\pp})^{\fv}$ is the space of invariants, and these are generated by the Chern forms of the Hodge bundles.}
\end{demo}

 The argument will  be given in several steps, some of which will carry over to the higher weight case.
\begin{itemize}\itemsep=0pt
\item[$(i)$] Since $D$
  is a symmetric space, the symmetry about $eV$ is given by an element in  $V$ which acts by $-1$ in the tangent space there.  It follows that the invariant forms are all of even degree.  In particular $\delta=0$, which gives the f\/irst statement in~\pref{4.9}.
  \item[$(ii)$] Next, using the unitary trick, and averaging, one infers that the space of invariant forms is isomorphic to the space of harmonic forms on the compact dual $\check D$, and this space is in turn isomorphic to $H^\ast(\check D)$.  This latter is known to be generated by the Chern classes of the universal vector bundle over $\check D$.
  \end{itemize}

  We will now show
  \begin{demo} \lab{4.10}
  {\bf Theorem.}
  {\it For a period domain, the invariant forms that are orthogonal to the {\bf algeb\-raic} ideal generated by $I$ and $\bar I$ are all of type $(p,p)$.}
  \end{demo}

Here, orthogonal means with respect to the canonical metrics constructed from the Hodge metrics in the Hodge bundles.  Working in the orthogonal space is equivalent to working modulo the ideal in the space of all forms.

 The f\/irst statement means that for a period domain the forms
 \[
 \lrb{\lrp{ \wedge^{\ast-1}\check{\pp}\wedge \scri_\C}^\bot}^\fv\]
 are all of type $(p,p)$.  This does not yet take into account the forms in $\wedge^{\ast-2}\check{\pp}\wedge\scri _2$.

 \begin{proof} At a reference point $F\in D$, we have a canonical representation
 \[
 \varphi:\ S^1\to \Aut(H_\R,Q)=G_\R \]
   obtained from a representation $\varphi_\C:\C^\ast\to\Aut(H_\C,Q)$ def\/ined by
 \[
 \varphi_\C(z)v = z^p\ol{z}^q v,\qquad  v\in H^{p,q} .\]
 Restricted to $S^1=\{|z|=1\}$, $\varphi_\C(z)=\varphi(z)$ is real.  We observe that $S^1$ is a subgroup of $V$; the Weil operator is just $C=\varphi(\sqrt{-1})$.

 On the Lie algebra $\scrG$, we have the induced Hodge structure
 \[
 \scrG_\C=\opplus_i\scrG^{-i,i} .\]
 We then have
 \[
 \fv_\C=\scrG^{0,0}\]
 and
 \[
 \pp_\C=\pp^+\oplus \pp^{-},\qquad  \pp^-=\ol{\pp}^+,\]
 where
 \[
 \left\{\begin{array}{l}
 \pp^+   = \opplus_{i>0}\scrG^{-i,i} \cong T_F^{(1,0)}D,\\[10pt]
 \pp^-   = \opplus_{i<0} \scrG^{-i,i}\cong T_F^{(0,1)}D.\end{array}\right.\]
 For the canonical Pfaf\/f\/ian system we have
 \[
 \scri _\C= \scri ^+ \oplus \scri ^-,\qquad \scri ^- = \ol{\scri }^+,\]
 where
 \[\left\{\begin{array}{l}
 \scri ^+   = \check{\scrG}^{-1,1}\subset \check{T}_F^{1,0}D,\\[4pt]
 \scri ^-   = \check{\scrG}^{1,-1}\subset \check{T}_F^{0,1}D .\end{array}\right.\]
  Thus the orthogonal complement to $\wedge^{\ast-1}\check{\pp}_\C
  \wedge \scri_\C$ are the forms of positive degree in
  \[
  \wedge^\ast \lrp{\check{\scrG}^{-1,1}\oplus \check{\scrG}^{1,-1}}\cong
  \wedge^\ast\check{\scrG}^{-1,1}\otimes \wedge^\ast \check{\scrG}^{1,-1} .\]
  On $\scrG^{-1,1}$, $\Ad\varphi(z)$ acts by $z^{-1} \ol{z}=z^2$; on $\scrG^{1,-1}$ it acts by $z\ol{z}^{-1}=z^{-2}$.  Thus on $\wedge^p\check{\scrG}^{-1,1}\otimes \wedge^q \check{\scrG}^{1,-1}$, $\Ad \varphi(z)$ acts by $z^{2(p-q)}$.  Since $\varphi(z)\in V$ for $z\in S^1$, we see that any $\ad\fv$-invariant forms must be of type $(p,p)$.
  \end{proof}

  As mentioned in the introduction, we expect that \pref{4.10} may be the f\/irst step in showing that the $G_\R$-invariant forms modulo the {\it differential} ideal generated by $I$ are generated by the Chern forms of the Hodge bundles\footnote{As mentioned in the introduction, this result has now been proved.  The integrability conditions play the central role in the argument, which is representation-theoretic in nature.}.  Since, except in the classical case when $I=0$, the invariant forms modulo the algebraic ideal generated by $I$ def\/initely contain more than the Chern forms, the integrability conditions will have to enter in an essential way.  In fact, they enter in the proof of the following
  \begin{demo}\lab{4.11}
 {\bf Theorem.} {\it The Chern forms of the Hodge filtration bundles satisfy}
 \[
 c_i(F^p)c_j(F^{n-p})=0\qquad \hensp{if} \quad i+j>h^{p,n-p} .\]
 \end{demo}

\begin{proof}  The argument will proceed in several  steps, the key one where integrability is used being~\pref{4.17} below.

 \medbreak\noindent {\it Step one.}  For a $d\times d$ matrix $A$ we set
 \[
 \chi_A(t) =:\det (A-tI)=t^d+c_1(A) t^{d-1}+\dots+c_d(A) .\]
 Then we have
 \begin{demo}\lab{4.12}
 {\bf Lemma.}  {\it If $A$, $B$ are $d\times d$ matrices with $AB=0$, then}
 \[
 c_i(A) c_j (B)=0\qquad \hbox{\it if }i+j>d .\]
 \end{demo}
 \begin{proof} We have
 \begin{gather*}
 \chi_A(t)\chi_B(u) =  \det(A-tI)\det (B-uI) = \det(-tB-uA+utI) ,\end{gather*}
 and all terms in $\det(-tB-uA+utI)$ have degree at least $d$ in $u$, $t$.  From
 \[
 \chi_A(t)\chi_B(u)=\sum_{i,j}c_i(A)c_j(B)t^{d-i}u^{d-j}\]
 we have $c_i(A)c_j(B)=0$ if $(d-i)+(d-j)<d$; i.e.\ if $i+j>d$.\end{proof}

 \medbreak\noindent {\it Step two.}
 We let $E\subset \scrG^{-1,1}$ be an integral element at a point $F\in D$.  Then we have $\check E$-valued maps
 \[
 F^p/F^{p+1}\mrar{A_p}
 F^{p-1}/F^p ,\]
 where   we think of $A_p$ as a matrix of size $h^{p,n-p}\times h^{p-1,n-p+1}$ whose entries are $\check E$.  The integrability conditions are
 \begin{equation}\lab{4.13}
 A_{p-1}\cdot A_p=0 ,\end{equation}
 where we are multiplying matrices using the wedge product of their entries.

 Using the dualities induced by $Q$ and denoting the transpose of a matrix $A$ by $A^\ast$ we have a~commutative diagram
 \[
 \begin{array}{ccc}
 \begin{picture}(10,15)
\put(-5,15){\line(3,-1){10}}
\put(5,11.75){\line(3,1){10}}
\end{picture}&&  \begin{picture}(10,15)
\put(-5,15){\line(3,-1){10}}
\put(5,11.75){\line(3,1){10}}
\end{picture}\\[-18pt]
(F^{p-1}/F^p)&\lrar{A^\ast_p}&(F^p/F^{p+1})\\
\wrvert&&\wrvert\\
F^{n-p+1}/F^{n-p+2}&\lrar{A_{n-p+1}}&F^{n-p}/F^{n-p+1} .\end{array}
\]

On $F^p$ we have an Hermitian metric induced by $(v,w)=Q(Cv,\bar w)$.

\begin{lem}
 Up to a constant, the curvature matrix of the metric connection is given by
 \begin{equation}\lab{4.14}
 \Theta_{F^p}=\begin{pmatrix}
 \ol{A}^\ast_pA_p&0\\[6pt]
 0&0\end{pmatrix} .\end{equation}\end{lem}

The notation means that we write the orthogonal direct sum decomposition
 \[
 F^p=H^{p,n-p}\oplus F^{p+1} .\]
 A consequence of \pref{4.14} is
 \begin{equation}\lab{4.15}
 c_k(F^p)=0\qquad \hensp{for} \quad k>h^{p,n-p} .\end{equation}

\begin{proof}[Proof of the lemma]  $F^p$ is a sub-bundle of the f\/lat bundle $F^0= H_\C$ with the induced metric.  In this situation it is known \cite{G} that the curvature of the metric connection in $F^p$ is up to a constant given
 by
 \[
 \overline{\Pi}^\ast_{F^p/H}\cdot \Pi_{F^p/H}, \]
 where $\Pi_{F^p/H}$ is the $2^{\rm nd}$ fundamental form of $F^p$ in $H_\C$.  In the case at hand the $2^{\rm nd}$ fundamental form may be identif\/ied with $A_p$.\end{proof}

 \begin{rem}
 There is one subtlety here.  Because of the sign alternation in the Hodge metrics
 \[
 (u,v)=(\sqrt{-1})^{p-q}Q(u,\ol{v})\qquad\qquad u,v\in H^{p,q}\]
 the usual principle that ``curvatures decrease on holomorphic sub-bundles'' does not hold for the Hodge bundles.  For example, $c_1(F^n)>0$ on $\Theta_{F^p}$.  However, the signs are not an issue for us here.\end{rem}

 \noindent {\it Step three.}  From \pref{4.13} we have
 \begin{equation} \lab{4.16}
 \ol{A}^\ast_p A_pA_{p+1}\ol{A}^\ast_{p+1}=0,
 \end{equation}
 where the multiplication of matrix entries is wedge product.  For notational simplicity we omit the blocks of zeroes in the $\Theta_{F^p}$'s, so that
 \begin{gather*}
 \Theta_{F^{p}} = \ol{A}^\ast_{p}A_{p}, \\
 \Theta^\ast_{F^{n-p}} =  -A_{p+1}\ol{A}^\ast_{p+1} .\end{gather*}
 Using \pref{4.16} this gives the remarkable consequence
 \begin{equation} \lab{4.17}
 \Theta_{F^p} \Theta^\ast_{F^{n-p}}=0\end{equation}
 of integrability.  Since $\chi_{\Theta^\ast_{F^{n-p}}}(t)=\chi_{\Theta_{F^{n-p}}}(t)$, by step one we see that $c_i(F^p)c_j(F^{n-p})=0$ if $i+j>h^{p,n-p}$.\end{proof}

 \section[``Expected'' dimension counts for integral manifolds of an EDS]{``Expected'' dimension counts for integral manifolds\\ of an EDS}\lab{5}

An important aspect in algebraic geometry is that of ``expected" dimension counts.  Informally and in f\/irst approximation, this means counting the number of parameters of solutions to a~system of algebraic equations, where ``expected" means assuming some sort of ``general position''.  When the solution varieties are also subject to dif\/ferential constraints, the problem changes character in an interesting way.  In this section we will discuss this   for the EDS arising from VHS's.

To frame the general issue we assume given a diagram of regular mappings of complex mani\-folds
\begin{equation}\lab{5.1}
\begin{array}{ccc}
X&\srar{f}&A\\
\cup&&\cup\\
Y&\srar{}&B\end{array}\end{equation}
where $f$ is an immersion and $Y=f^{-1}(B)$.  Then
\begin{gather*}
\hbox{``expected''} \ \codim_XY = \codim_A B =\mathrm{rank}(TA/TB) ,\end{gather*}
where it is understood that $TA/TB$ is restricted to $B$.
If $Y$ is non-empty, then the actual codimension satisf\/ies
\begin{equation}\lab{5.2}
\codim_X Y\leqq \mathrm{rank}(TA/TB) ,\end{equation}
with equality holding when $f(X)$ meets $B$ transversely.

Now suppose that $I\subset T^\ast A$ is a holomorphic Pfaf\/f\/ian system and  $f:X\to A$ is an integral manifold of $I$.  Let $W=I^\bot$ be the corresponding distribution.  Then for the normal bundles we have that
\[
f_\ast : \ TX/TY \to W/W\cap TB\]
is injective, so that the above may be ref\/ined to
\begin{equation}\lab{5.3}
\codim_X Y\leqq {\rm rank}(W/W\cap TB)\leqq {\rm rank}(TA/TB) .\end{equation}
Informally we may say that: {\it Subjecting $f:X\to A$ to a differential constraint {\bf decreases} the codimension of $Y=f^{-1}(A)$ in $X$.}  By ``decreases" we mean that $\codim_XY$ is less than the ``expected" codimension in  the absence of dif\/ferential constraints, as explained above.

However, {\it when we take into account integrability a still further refinement of $\pref{5.3}$ may be expected.}  This is because in general integral elements of $I$ may be expected to have dimension smaller, frequently much less, than rank $W$.

Rather than discuss the general aspects of this, we turn to a special case that is geometrically motivated.  Let $D$ be a period domain for polarized Hodge structures of even weight $n=2m$.  At a reference Hodge structure $F\in D$ we let $\zeta \in H_\R\cap H^{m,m}$ be a real vector of type $(m,m)$.

\begin{defi}
{\it The  Noether--Lefschetz locus is the subvariety $D_\zeta\subset D$ where $\zeta$ remains of type $(m,m)$.}
\end{defi}

\begin{demo}\lab{5.4}
{\bf Proposition.} {\it Let $G_\zeta\subset G_\R$ be the subgroup fixing $\zeta$ up to scaling.  Then $D_\zeta=G_\zeta \cdot F$ is the $G_\zeta$-orbit of $F$.  It is a homogeneous complex sub-manifold of $D$ of codimension given by}
\[
\codim_D D_\zeta=h^{(2m,0)}+\cdots +h^{ (m+1,m-1)} .\]
\end{demo}

\begin{proof} This is a matter of routine checking.  Setting $h_k=h^{(2m-k,k)}$, we have f\/irst
\[
\left\{\begin{array}{l}
G_\R   \cong O(a,b), \qquad a+b=h_0+\dots+h_{2m},\\
V   \cong \scrU(h_0)\times\dots\times \scrU(h_{m-1})\times O(h_m) . \end{array}\right.\]
Next, depending on whether $m$ is even or odd, we have $G_\zeta\cong O(a-1,b)$ or $O(a,b-1)$.  Finally, the same linear algebra argument that shows that $G_\R$ acts transitively on $D$ shows that $G_\zeta$ acts transitively on $D_\zeta$ and
\[
D_\zeta = G_\zeta/V_\zeta,
\]
where
\[
V_\zeta=G_\zeta\cap V\cong \scrU(h_0)\times\dots\times \scrU(h_{m-1})\times O(h_{m}-1) .\]
From this we may conclude the above codimension count.
\end{proof}

For f\/ixed $\zeta\in H_\R\cap H^{(m,m)}$ and at a variable point in $D$ we write the Hodge decomposition of $\zeta$ as
\[
\left\{ \begin{array}{l}
\zeta = \zeta^{2m,0} +\cdots+ \zeta^{m,m}+\dots+\zeta^{0,2m},\\[4pt]
\zeta^{\overline{2m-p,p}}=\zeta^{p,2m-p} .\end{array}\right.\]
Then $D_\zeta$ is def\/ined by the equations
\[
\zeta^{m-1,m+1}=\dots= \zeta^{0,2m}=0 .\]
The above proposition says that these  equations are independent and def\/ine $D_\zeta$ as a smooth complex  submanifold of $D$.

At a point $F\in D$ we let
\[
\left\{\begin{array}{l}
E   \subset T_F D\quad\hbox{be an integral element of $I$},\\[4pt]
E_\zeta   = E\cap T_FD_\zeta .\end{array}\right.\]
For $\varphi\in \scrG^{-1,1}\subset \oplus \Hom(H^{2m-p,p},H^{2m-p-1,p+1})$ we write $\varphi=\varphi_0+\dots+\varphi_{2m-1}$ where $\varphi_p\in \Hom(H^{2m-p,p},H^{2m-p-1,p+1})$ and  $\varphi_p$ and $\varphi_{2m-p-1}$ are dual.  Then
\[
E_\zeta=\{\varphi\in E:\varphi_m (\zeta)=0\hensp{in} H^{m-1,m+1}\} .\]
This is equivalent to
\begin{equation}\lab{5.5}
E_\zeta=\{\varphi\in E:Q(\eta,\varphi(\zeta))=0\hensp{for all}\eta\in H^{m+1,m-1}\} .\end{equation}
Thus, without taking the integrability conditions into account we have
\[
\codim_E E_\zeta\leqq h^{m+1,m-1} .\]

However, due to the integrability conditions the equations \pref{5.5} may not be independent.  In order to illustrate the essential point, we begin by considering the f\/irst non-trivial case $m=2$.  For any $\varphi\in E_\zeta$, $\psi\in E$ and $\om \in  H^{4,0}$, using that $E$ is an integral element so that $\varphi$ and $\psi$ commute,
\begin{gather*}
Q(\varphi(\om),\psi(\zeta)) =  -Q(\psi\varphi(\om),\zeta)
 =  -Q(\varphi\psi(\om),\zeta)
 =  Q(\psi(\om),\varphi(\zeta))
 = 0  ; \end{gather*}
i.e., for each $\psi\in E$ the linear equations
\[
Q(\eta,\psi(\zeta))=0,  \qquad \eta\in H^{3,1}\]
that def\/ine the condition that $\psi\in E_\zeta$ are decreased in rank by
\begin{equation}\lab{5.6}
\sigma_\zeta=: \dim \rim \{E_\zeta\otimes H^{4,0}\to H^{3,1}\}  .\end{equation}
We thus have the
\begin{demo}\lab{5.6p}
{\bf Proposition.} {\it For $\sigma_\zeta$ defined as above}
\[
\codim_E E_\zeta \leqq h^{1,3}-\sigma_\zeta .\]
\end{demo}

The general case goes as follows:  With $E_\zeta$ def\/ined as above, for each $p$ with $0\leqq p\leqq m-2$ we consider the maps
\[
\kappa^p_\zeta:E \otimes \Sym^{p-m+1}E\otimes H^{2m-p,p}\to H^{m+1,m-1}\]
and we set
\begin{equation}\lab{5.8}
\sigma_\zeta=\dim\lrc{ \hbox{span(Images }\kappa^p_\zeta)\hensp{for} 0\leqq p\leqq m-2} .
\end{equation}
Then the straightforward extension of the above argument gives the
\begin{demo}\lab{5.9}
{\bf Proposition.} {\it For $\sigma_\zeta$ defined by $\pref{5.8}$, we have}
\[
\codim_E E_\zeta\leqq h^{m-1,m+1}-\sigma_\zeta .\]
\end{demo}

In algebro-geometric terms, this says that {\it the ``richer" the multiplicative structure in the $1^{\rm st}$ order variation of the Hodge structure, the smaller the codimension of Noether--Lefschetz loci.}

\begin{rem}
The above is predicated on the implicit assumption that, for a variation of Hodge structure $S\to \Gamma\bsl D$, we have that the Noether--Lefschetz locus $S_\zeta\subset S$ is reduced, so that for   general points of a component of $S_\zeta$ we have $\codim_S(S_\zeta)=\codim_{T_s S}(T_sS_\zeta)$.  For far as we know, there are not yet any examples coming from algebraic geometry where this assumption is not satisf\/ied, although we feel that such examples may be expected.
\end{rem}

\begin{exam}  We consider a smooth hypersurface
\[
X\subset\P^5\]
of degree $d\geqq 6$ and which contains a 2-plane $P\cong\P^2$.
\end{exam}

\begin{demo} \lab{5.10} {\bf Theorem.} {\it For the primitive part of the fundamental class $ [P]\in H^4(X,\Z)$ of $P$, we have equality in~$\pref{5.9}$.}
\end{demo}

\begin{proof} We will denote by $V\cong\C^6$ a vector space with coordinates $x_1,\dots,x_6$ such that $X\subset \P\check{V}$ is given by an equation
\[
F(x)=0\]
where $F\in V^d=: \Sym^d V$ is a homogeneous polynomial of degree $d\geqq 6$.  It is well-known, and will be proved below, that
\begin{demo} \lab{5.11}
{\it  At a general $X$ containing a $2$-plane $P$, the locus of all degree $d$ hypersurfaces $X'$ close to~$X$ and containing a $2$-plane $P'$ close to $P$ is smooth and of codimension $\frac{(d+1)(d+2)}{2}- \dim \Gr(3,6)$ in the space of all degree $d$ hypersurfaces in $\P \check V$.}
\end{demo}

Let $P$ be given by $x_1=x_2=x_3=0$ so that
\[
F(x) = \sum^3_{i=1}x_iG_i(x),
\]
where $G_i\in V^{d-1}$.  Denote by
\[
F_t(x)=F(x)+t \dot{F}(x),\qquad  \dot{F}\in V^d, \]
a $1^{\rm st}$ order perturbation of $F$.  The condition that $P$ move to $1^{\rm st}$ order to a
2-plane $P_t\subset X_t=\{F_t(x)=0\}$ is that
\[\dot{F} = \sum_i x_i\dot{G}_i+\sum_il_iG_i,
\]
where $\dot{G}_i\in V^{d-1}$ is the tangent to a $1^{\rm st}$ order variation of $G_i$ and the $l_i\in V$ are linear forms.  We will use the notation $(X',P')$ for  the $1^{\rm st}$ order perturbation of $(X ,P )$.

We will denote by $(H_1{,}{\dots}{,}H_m)$ the ideal generated by forms $H_1{,}{\dots}{,}H_m$, and  by $(H_1{,}{\dots}{,}H_m)_k$ the degree $k$  part of that ideal.  Thus, the subspace of $V^d$ that gives the $1^{\rm st}$ order deformations $(X',P')$ of $(X,P)$ is
\[
(x_1,x_2,x_3,G_1,G_2,G_3)_d .\]

The {\it Jacobian ideal} is
\[
J=(\part_{x_1}F,\dots,\part_{x_6}F) .\]
For references to Jacobian ideals and the polynomial description of the cohomology of hypersurfaces we suggest \cite{A,CD}, and \cite{D}.
  The tangent space to $1^{\rm st}$ order deformations of the projective equivalence class of $X$ is given by
\[
T=V^d/J_d .\]
With this identif\/ication, by what was said above the subspace $T_P\subset T$ giving $1^{\rm st}$ order
deformations of equivalence classes of pairs $(X,P)$ is given by
\begin{equation}
\lab{5.12}
T_P=(x_1,x_2,x_3,G_1,G_2,G_3)_d /J_d .\end{equation}

A basic identif\/ication is
\begin{equation}\lab{5.13}
H^{4-p,p}(X)_\prim = V^{(p+1)d-6}/J_{(p+1)d-6} .\end{equation}
Thus, for example
\[
H^{4,0}(X)\cong V^{d-6} .\]
Using the identif\/ication $T=V^d/J_d$, the dif\/ferential of the period mapping
\[
T\to \opplus_p\Hom(H^{4-p,p}(X),H^{4-p-1,p+1}(X))\]
is given by multiplication in the ring $V^\bullet/J_\bullet$; i.e.\ by
\begin{equation}  \lab{5.14}
V^d/ J_d\otimes V^{(p+1)d-6}/J_{(p+1)d-6}\to V^{(p+2)d-6}/J_{(p+2)d-6} .\end{equation}

We denote by $\zeta\in H^{2,2}(X)_{\rm prim}$ the primitive part of the fundamental class $[P]$ of $P$.

 \begin{demo}\lab{5.15}
 {\bf Proposition.}  {\it If we define
 \[
 T_\zeta=\{ H\in V^d/J_d: H\cdot\zeta=0\hensp{in}V^{4d-6}/J_{4d-6}\}
 \]
 then}
 \[
 T_\zeta=T_P .\]
 \end{demo}

 The inclusion $T_P\subseteq T_\zeta$ is clear geometrically: It means that if $(X,P)$ deforms to $1^{\rm st}$ order to $(X',P')$ as above, then $\zeta$ deforms to a Hodge class $\zeta'\in H^2(X')_{\prim}$.  In fact, $\zeta' = [P']_{\prim}$.  The proof of the reverse inclusion will come out indirectly from the argument to be given below.  The assertion \pref{5.11} is a consequence of \pref{5.15}.

 We set
 \[
 V_P=V\mid_P=V/(x_1,x_2,x_3) .\]
 The proof of the theorem and the proposition will be based on known commutative algebra properties of the  rings $V^\bullet/J_\bullet$ and $V^\bullet_P/J_{P,\bullet}$; namely there are perfect pairings
 \begin{equation}\lab{5.16}
 \left\{
 \begin{array}{ll}
(i)  & V^k/J_k\otimes V^{6d-12-k}/J_{6d-12-k}\to V^{6d-12}/J_{6d-12}\cong\C,\\[4pt]
(ii)&V^k_P/(G_1,G_2,G_3)_{P,k}\otimes V_P^{3d-6-k}/(G_1,G_2,G_3)_{P,3d-6-k}\\[4pt]
 &  \quad \to V^{3d-6}_P/(G_1,G_2,G_3)_{P,3d-6}\cong\C .\end{array}\right.\end{equation}
 Here, we denote by $(G_1,G_2,G_3)_{P,k}$ the degree $k$  part of the ideal generated by the $G_i\vert P$. The reason for \pref{5.16} is that $(x_1,x_2,x_3,G_1,G_2,G_3)$ is a {\it regular sequence} on $\P^5$, and this then implies that $G_1\vert_P, G_2\vert_P,G_3\vert_P$ is a regular sequence on $P$.  Then in general if $f_1,\dots,f_n$ is a regular sequence on $\P^n$ where $\deg f_i=d_i$, there is a perfect pairing $V^a/(f_1,\dots,f_n)_a  \otimes V^b/(f_1,\dots,f_n)_b\to V^{a+b}/(f_1,\dots,f_n)$ where $a+b=\sum_i d_i-n$.

 The argument will proceed in f\/ive steps, the f\/irst of which is \pref{5.13} above.

 \medbreak\noindent {\it Step two.}  We have
 \[
 T_P H^{4,0}=(x_1,x_2,x_3,G_1,G_2,G_3)_{2d-6}/J_{2d-6}\subseteq V^{2d-6}/J_{2d-6} .\]
 This follows from \pref{5.13} and \pref{5.14} in the case $p=0$.

 \medbreak\noindent {\it Step three.} Denoting by $[P]$ the fundamental class of $P$, the map
 \[
 H^{2,2}(X)\srar{[P]} H^{4,4}(X)\cong\C\]
 may, using $(ii)$ in \pref{5.14} which gives an isomorphism $V^{3d-6}/(G_1,G_2,G_3)_{P,3d-6}  \cong\C$, be identif\/ied with
 \[
 V^{3d-6}/J_{3d-6}\to V^{3d-6}_P/ (G_1,G_2,G_3)_{P,3d-6}\cong\C .\]

\begin{proof}  Since $X$ is non-singular the restrictions $G_i\vert_P$ have no common zeroes and hence form a~regular sequence.  The corresponding Koszul resolution of $\scrO_P$ then gives
 \[
 0\to \scrO_P(-3(d-1))\to\opplus_3 \scrO_P(-2(d-1))\to\opplus_3\scrO_P(-(d-1))\to\scrO_P\to 0 .
 \]
 Tensoring with $\scrO_P(3d-6)$ gives, using $\Om^2_P\cong\scrO_P(-3)$
 \[
 0\to\Om^2_P\to\opplus_3\scrO_P(d-4)\to\opplus_3 \scrO_P(2d-5)\to\scrO_P(3d-6)\to 0 ,\]
 from which we infer that
 \[
 V^{3d-6}_P/(G_1,G_2,G_3)_{P,3d-6} \cong H^2(P,\Om^2_P)\cong\C .\]
 Moreover, under this isomorphism a generator of $V^{3d-6}_P/(G_1,G_2,G_3)_{P,3d-6}$ maps to the fundamental class.  Since the map $H^{2,2}(X)\srar{[P]} H^{4,4}(X)$ is given by evaluating a class in $H^{2,2}(X)$ on the fundamental class of $P$, by standard arguments we may infer the assertion in Step~3.
\end{proof}

 \medbreak\noindent {\it Step four.} We f\/irst observe that the map
 \begin{equation}\lab{5.17}
  \zeta \otimes V^d/J_d\otimes T_P\otimes V^{d-6}\to V^{4d-6}/J_{4d-6}\cong\C\end{equation}
 is zero.  Here, we recall that    $\zeta\in V^{3d-6}/J_{3d-6}$ is the primitive part of $[P] \in H^{2,2}(X)$.  We are using  \pref{5.14} that the action on cohomology of tangent vectors to deformations of equivalence classes of hypersurfaces in $\P^5$  is given by multiplication in the ring $V^\bullet/J_\bullet$.   The fact that the above map is zero results from the def\/inition of $T_\zeta$ as the kernel of the map
 \[
 \begin{array}{ccc}
 T&\mrar{ \zeta }&H^{1,3}\\[5pt]
 \wrvert&&\wrvert\\[5pt]
 V^d/J_d&\mrar{}&V^{4d-6}/J_{4d-6} .\end{array}\]
 Then the claim is that the map \pref{5.17} may be identif\/ied with
 \begin{equation}\lab{5.18}
 T_P \cdot V^{d-6}\mid_P\in V^{3d-6}_P/(G_1,G_2,G_3)_{P,3d-6} \cong\C .\end{equation}
 This follows from Step three above.

 \medbreak\noindent {\it Step five.} We now put everything together.  To prove the theorem it will suf\/f\/ice to show that $T_P$ is the kernel of
 \[
 V^d/J_d \srar{ \zeta } V^{4d-6}/J_{4d-6} .\]
 This is because for all $S\in V^d/J_d$,
 \begin{equation}
 \lab{5.19}
 R\cdot\zeta=0\hensp{in} V^{6d-6}/J_{6d-6}\hensp{for all}R\in V^{2d-6}\end{equation}
 is, by \pref{5.16}$(i)$, equivalent to
 \[
 S\cdot  \zeta =0\hensp{in}V^{4d-6}/J_{4d-6} .\]
 By step four, \pref{5.19} is in turn equivalent to
 \[
 RS\mid_P \in (G_1,G_2,G_3)_{P,3d-6} ,\]
 which by \pref{5.16}$(ii)$ is the same as the condition
 \[
 S\mid_P\in (G_1,G_2,G_3)_{P,d} .\]
 This last statement is easily seen to be equivalent to
\begin{gather*}
S\in T_P.\tag*{\qed}
\end{gather*}\renewcommand{\qed}{}
\end{proof}

 \subsection*{Reprise} We consider the case that arises in the case of a family of Calabi--Yau fourfolds.  Thus, we assume that $h^{4,0}=1$ and denote by $\om\in H^{4,0}$ a generator giving an isomorphism $H^{4,0}\cong \C$.  Moreover, let $T\subset \scrG^{-1,1}$ be an integral element and, as would be the case for Calabi--Yau's, we assume that the map  given by \pref{5.13}
 \[
 T\to \Hom(H^{4,0},H^{3,1})\cong H^{3,1}\]
 is an isomorphism.  Using this map we may identify $T$ with $H^{3,1}$ and denote its dimension by~$m$ (for moduli).

 Each $\zeta\in H^{2,2}$ def\/ines a quadric
 \[
 Q_\zeta\in \Sym^2\check T\]
 given for $\theta,\theta'\in T$ by
 \[
 \lra{Q_\zeta,\theta\cdot\theta'}=Q(\theta\cdot \theta'\om,\zeta).\]
 The fact that this is symmetric in $\theta$ and $\theta'$ is because $T$ is an integral element.

 Denote by $T_\zeta\subset T$ the intersection of $T$ with the tangent space to the Noether--Lefschetz locus $D_\zeta\subset D$.

 \begin{demo}\lab{5.20}
 {\bf Proposition.} {\it We have
 \[
 \codim T_\zeta=\rank Q_\zeta .\]
 Moreover,  in the case of the family of hypersurfaces  $X \subset \P^5$ of degree six and where $\zeta$ is the primitive part of the fundamental class of a plane $P\subset X$, we have $T_P=T_\zeta$.}
 \end{demo}

 \begin{proof}
 We view $Q_\zeta$ as a map
 \[
 Q_\zeta: \ T\to\check T .\]
 Then the proof of Proposition~\ref{5.6p} gives
 \[
 T_\zeta= (\rim Q_\zeta)^\bot ;\]
 thus
 \[
 \codim T_\zeta=\rank Q_\zeta\]
 as desired. \end{proof}

 For the case of $P\subset X\subset \P^5$, $\deg X=6$, we have
 \[
 T=V^6/J_6 .\]
 Moreover, we have seen in step three in the proof of \pref{5.15} that
 \[
 T_\zeta=\ker\lrc{ V^6/J_6\to V^6_P/(G_1,G_2,G_3)_{P,6}} .
\]
Thus
\[
\rank Q_\zeta=\dim V^6_P/(G_1,G_2,G_3)_{P,6} .\]
Now, from the Koszul calculation above
\[
0\to \opplus_3 V_P \lrar{(G_1,G_2,G_3)} V^6_P
\]
is exact.  Thus
\begin{gather*}
\rank Q_\zeta =  \dim V^6_P-3\dim V_P = {8\choose 2} - 3\cdot 6  = 19.\end{gather*}
On the other hand, the number of conditions for $X$ to contain a $2$-plane is
\[
\dim V^6_P-\dim \Gr(3,6) = {8\choose 2}-9=19 .
\]

\begin{rem}
Let $X$ be a Calabi--Yau fourfold and $\zeta\in \Hg^2(X)_\prim$ a Hodge class.  We then have
\begin{demo}
\lab{5.21}
{\it If the Hodge conjecture is true and $\rank Q_\zeta=h^{3,1}$ is maximal, then $X$ is  defined over a~number field.}
\end{demo}

This is because if $Q_\zeta$ is non-singular, then the Noether--Lefschetz locus for $\zeta$ will be  $0$-dimensional.  If $\zeta=[Z]$ is the class of an algebraic cycle, then   $X$ is def\/ined over a f\/ield $k$ of transcendence degree $\geqq 1$, and by standard arguments we may, after passing to a f\/inite f\/ield extension, assume that $Z$ is also def\/ined over $k$.  The {\it spread} of $(X,Z)$ will then give a positive dimensional component to the Noether--Lefschetz locus of $\zeta$.

To  disprove the consequence of the Hodge conjecture that {\it Hodge classes are absolute}, it would be  suf\/f\/icient to f\/ind a Calabi--Yau fourfold not def\/ined over a number f\/ield and a Hodge class~$\zeta$ such that $Q_\zeta$ is non-singular.   \end{rem}

\addcontentsline{toc}{section}{References}
\LastPageEnding


\begin{thebibliography}{99}

\footnotesize\itemsep=0pt

 \bibitem{A}
 Allaud E.,
 Nongenericity of variations of Hodge structure for hypersurfaces of high degree,
 {\it Duke Math.\ J.} {\bf 129} (2005), 201--217,
\href{http://arxiv.org/abs/math.AG/0503346}{math.AG/0503346}.

 \bibitem{BCGGG}
 Bryant R.L., Chern S.S., Gardner R.B., Goldschmidt H.L., Grif\/f\/iths P.A.,
 Exterior dif\/ferential systems, {\it Mathematical Sciences Research Institute Publications}, Vol.~18, Springer-Verlag, New York, 1991.

 \bibitem{B-G}
 Bryant  R.L., Grif\/f\/iths P.A.,
 Some observations on the inf\/initesimal period relations for regular threefolds with trivial canonical bundle,
in  Arithmetic and Geometry, Vol.~II, {\it Progr.\ Math.}, Vol.~36, Birkh\"auser Boston, Boston, MA, 1983, 77--102.

 \bibitem{BG1}
 Bryant R.L, Grif\/f\/iths P.A.,
 Characteristic cohomology of dif\/ferential systems. I.~General theory,
 {\it J. Amer.\ Math.\ Soc.} {\bf 8} (1995), 507--596.

 \bibitem{BG2}
 Bryant R.L, Grif\/f\/iths P.A.,
 Characteristic cohomology of dif\/ferential systems. II.~Conservation law for a~class of parabolic equations,
 {\it Duke Math.\ J.} {\bf 78} (1995), 531--676.

 \bibitem{C}
 Carlson J.A., Bounds on the dimension of variations of Hodge structure,
 {\it Trans.\ Amer.\ Math.\ Soc.} {\bf 294} (1986), 45--64,
 Erratum, {\it Trans.\ Amer.\ Math.\ Soc.} {\bf 299} (1987), 429.

 \bibitem{CD}
 Carlson  J.A., Donagi R.,
 Hypersurface variations are maximal.~I,
 {\it Invent.\ Math.} {\bf 89} (1987), 371--374.

 \bibitem{CKT}
 Carlson J.A., Kasparian A., Toledo D.,
 Variations of Hodge structure of maximal dimension,
 {\it Duke Math.~J.} {\bf 58} (1989), 669--694.

 \bibitem{CS}
Carlson  J.A., Simpson C.,
Shimura varieties of weight two Hodge structures,
in Hodge Theory (Sant Cugat, 1985),
{\it Lecture Notes in Math.}, Vol.~1246, Springer, Berlin, 1987, 1--15.

\bibitem{CT}
 Carlson J.A., Toledo D.,
 Generic integral manifolds for weight-two period domains,
 {\it Trans.\ Amer.\ Math.\ Soc.} {\bf 356} (2004), 2241--2249,
 \href{http://arxiv.org/abs/math.AG/0501078}{math.AG/0501078}.

\bibitem{CT2}
Carlson  J.A., Toledo D.,
Variations of Hodge structure, Legendre submanifolds and accessibility,
{\it Trans.\ Amer.\ Math.\ Soc.} {\bf 311} (1989), 391--411.


 \bibitem{C-MS-P}
 Carlson  J.A., M\"uller-Stach S., Peters C.,
  Period mappings and period domains, {\it Cambridge Studies in Advanced Mathematics}, Vol.~85,  Cambridge University Press, Cambridge, 2003.

 \bibitem{D}
 Donagi R.,
 Generic Torelli for projective hypersurfaces,
 {\it Compositio Math. } {\bf 50} (1983), 325--353.

 \bibitem{Gr}
 Green M.,
 Koszul cohomology and the geometry of projective varieties.~II,
 {\it J. Differential Geom.} {\bf 20} (1984), 279--289.


   \bibitem{GG2}
  Green   M., Grif\/f\/iths P.,
  Algebraic cycles and singularities of normal functions. II, in Inspired by S.S.~Chern,
  {\it Nankai Tracts Math.}, Vol.~11, World Sci. Publ., Hackensack, NJ, 2006, 2006, 179--268.

 \bibitem{G}
 Grif\/f\/iths P.,
 Hermitian dif\/ferential geometry and the theory of positive and ample holomorphic vector bundles,
 {\it J. Math.\ Mech.} {\bf 14} (1965), 117--140.


    \bibitem{IL}
   Ivey  T.A., Landsberg J.M.,
Cartan for beginners: dif\/ferential geometry via moving frames and exterior dif\/ferential systems, {\it Graduate Studies in Mathematics}, Vol.~61, American Mathematical Society, Providence, RI, 2003.

    \bibitem{M}
Mayer R.,
Coupled contact systems and rigidity of maximal dimensional variations of Hodge structure,
{\it Trans.\ Amer.\ Math.\ Soc.} {\bf 352} (2000), 2121--2144,
\href{http://arxiv.org/abs/alg-geom/9712001}{alg-geom/9712001}.

    \bibitem{O1}
Otwinowska A.,
Composantes de petite codimension du lieu de Noether--Lefschetz: un argument asymptotique en faveur de la conjecture de Hodge pour les hypersurfaces,
{\it J. Algebraic Geom.} {\bf 12} (2003), 307--320.

    \bibitem{O2}
Otwinowska A.,
Composantes de dimension maximale d'un analogue du lieu de Noether--Lefschetz,
{\it Compositio Math.} {\bf 131} (2002), 31--50.

    \bibitem{V}
Voisin C.,
Hodge loci and absolute Hodge classes,
{\it Compositio Math.} {\bf 143} (2007), 945--958,
\href{http://arxiv.org/abs/math.AG/0605766}{math.AG/0605766}.

  \end{thebibliography}
  \end{document}